\documentclass[a4paper,12pt]{article}

\usepackage{amsmath,amsfonts,amssymb}
\usepackage[colorlinks]{hyperref}

\newtheorem{theorem}{Theorem}[section]
\newtheorem{lemma}[theorem]{Lemma}

\newtheorem{proposition}[theorem]{Proposition}

\newtheorem{remark}[theorem]{Remark}
\newtheorem{example}[theorem]{Example}
\newtheorem{definition}[theorem]{Definition}

\numberwithin{equation}{section}

\allowdisplaybreaks

\begin{document}

\title{%\vspace*{-1.5cm}
Approximations  of Stochastic Navier-Stokes Equations}

\author{ Shijie Shang$^{1}$ \and
Tusheng Zhang$^{2}$}

\footnotetext[1]{\ School of Mathematical Sciences, University of Science and Technology of China, 230026 Hefei, China. Email: ssjln@mail.ustc.edu.cn}
\footnotetext[2]{\ School
of Mathematics, University of Manchester, Oxford Road, Manchester
M13 9PL, England, U.K. Email: tusheng.zhang@manchester.ac.uk}

\maketitle

\begin{abstract}
In this paper we show that solutions of two-dimensional stochastic Navier-Stokes equations driven by Brownian motion can be approximated by stochastic Navier-Stokes equations forced by pure jump noise/random kicks.
\end{abstract}


\section{Introduction}

Stochastic Navier-Stokes equations (SNSEs) are now a widely accepted model for fluid motion with random perturbations.
In this paper, we consider the two-dimensional stochastic Navier-Stokes equations with Dirichlet boundary conditions on a bounded domain, which is given as follows:
\begin{align}\label{1.a}
\left\{
\begin{aligned}
& du -\kappa \Delta u\,dt+(u\cdot\nabla)u\,dt+\nabla\mathfrak{P}\,dt = F(u)\,dt+\sum_{i=1}^{m}\sigma^{i}(u)\,dW^{i},\quad \rm{ in }\ \mathcal{O}\times(0,T], \\
&\begin{aligned}
& {\rm{div}}\,u=0 \quad &&\rm{in}\ \mathcal{O}\times(0,T], \\
& u=0  &&\rm{in}\ \partial \mathcal{O}\times(0,T], \\
& u(0)=h  &&\rm{in}\ \mathcal{O},
&\end{aligned}
\end{aligned}
\right.
\end{align}

\noindent where $\mathcal{O}$ is a bounded domain of $\mathbb{R}^2$ with boundary $\partial \mathcal{O}$ of class $\mathcal{C}^{3}$. $u=(u_1,u_2)$ and $\mathfrak{P}$ represent the random velocity and modified pressure, respectively. $\kappa$ is the kinematic viscosity, for simplicity, we let $\kappa=1$ in this paper. $W=(W^{1}(t),\cdots,W^{m}(t))$ is a $m$-dimensional standard Brownian motion. The fluid is driven by external force $F(u)\,dt$ and the random noise $\sum_{i=1}^{m}\sigma^{i}(u)\,dW^{i}$.

\vskip 0.3cm

%Stochastic Navier-Stokes equations has been initiated by Bensoussan and Temam in \cite{1973-Bensoussan-p195-222}.

Stochastic Navier-Stokes equations have been studied by many people. There is a great amount of literature. Let us mention a few. SNSEs driven by white noise in time were first studied by Bensoussan and Temam in \cite{1973-Bensoussan-p195-222}.
%Stochastic Navier-Stokes equation have been intensively studied since the work of Bensoussan and
%Temam in [???].driven by white noise type random force
%For 2-D SNSEs on a bounded domain, there are many literature.
%has also been studied by many authors.
The existence and uniqueness of solutions of 2-D SNSEs driven by L\'{e}vy noise were obtained in \cite{2014-Brzezniak-p283-310}, large derivation and moderate deviation principles were established in \cite{2015-Zhai-p2351-2392,2017-Dong-p227-254}.
The ergodic properties and invariant measures of the 2-D SNSEs were studied in \cite{2006-Hairer-p993-1032} and \cite{1994-Flandoli-p403-423}.

\vskip 0.3cm

The aim of this paper is to study the approximations of  SNSEs in (\ref{1.a}) by SNSEs forced by Poisson random measures. One of the motivations is to shine some light on numerical simulations of SNSEs driven by pure jump noise. Recently, Nunno and Zhang in \cite{2016-DiNunno-p1443-1466} obtained such an  approximation for a general class of SPDEs. However, the results in \cite{2016-DiNunno-p1443-1466} couldn't cover the stochastic Navier-Stokes equations, an important model in fluid dynamics. The difficulty lies in establishing the tightness of the approximating equations in the space of Hilbert space-valued right continuous paths with left limits. To overcome this difficulty, we first assume that the initial value has higher regularity, the external force and the coefficients of the jump noise take values in a more regular space, so that we can derive an uniform estimate of the stronger norm of the approximating solutions. With these estimates, we are able to prove the tightness of the approximating equations by Aldou's criterion, then through martingale characterization we show that the limit of the solutions of approximating equations is the solution of the SNSE driven by Brownian motion. We emphasize that the  method of establishing the tightness here is different and simpler than that used in \cite{2016-DiNunno-p1443-1466}.  In the second step, we are able to remove the regularity restrictions on the coefficients and the initial condition by using finite dimensional approximations and establishing some uniform convergence in probability of the approximating solutions. In the final part of the paper, we provide several illustrating examples.
%Deriving the special coefficients used in \cite{2016-DiNunno-p1443-1466} from the class of examples shows that under weaker conditions we establish our results.
%To the best of our knowledge, there is only one paper this is the second paper to consider such approximations for SPDEs.

%studied the approximation problem under the general framework of SPDEs driven by a special pure jump noise, the results in there , while, in our paper, we solve the SNSE driven by general jump coefficients under weaker conditions.
%To the best of our knowledge

\vskip 0.3cm
The rest of the paper is organized as follows. In Section 2 we lay down the precise framework. The main part is Section 3, where the approximations are established. In Section 4  some examples are provided.

%%%%%%%%%%%%%%%%%%%%%%%%%%%%%%%%%%%%%%%%%%%%%%%%%%%%%%%%%%%%%%%%%%%%%%%%%%%%%%%%%%%%%%%%%%%%%%%%%%%%%%%%%%%%%%%
%%%%%%%%%%%%%%%%%%%%%%%%%%%%%%%%%%%%%%%%%%%%%%%%%%%%%%%%%%%%%%%%%%%%%%%%%%%%%%%%%%%%%%%%%%%%%%%%%%%%%%%%%%%%%%%
%%%%%%%%%%%%%%%%%%%%%%%%%%%%%%%%%%%%%%%%%%%%%%%%%%%%%%%%%%%%%%%%%%%%%%%%%%%%%%%%%%%%%%%%%%%%%%%%%%%%%%%%%%%%%%%
%%%%%%%%%%%%%%%%%%%%%%%%%%%%%%%%%%%%%%%%%%%%%%%%%%%%%%%%%%%%%%%%%%%%%%%%%%%%%%%%%%%%%%%%%%%%%%%%%%%%%%%%%%%%%%%

\section{Framework}

\noindent Let $(\Omega, {\cal F}, P)$ be a probability
space equipped with a filtration $\{{\cal F}_t\}_{t\geq 0}$ satisfying the usual conditions.
%Let $W(t)$ be a one-dimensional ${\cal F}_t$-Brownian motion.
$\nu^{i} (dx), i=1,\cdots,m$ denote $\sigma$-finite measures on the measurable space $(\mathbb{R}_0, {\cal B}(\mathbb{R}_0))$, where $\mathbb{R}_0:=\mathbb{R}\setminus \{0\}$.
Let $N^{i}, i=1,\cdots,m$ be mutually independent ${\cal F}_t$-Poisson random measures on $[0,T]\times \mathbb{R}_0$ with intensity measure $dt\times\nu^{i}(dz)$ respectively.
For $U\in {\cal B}(\mathbb{R}_0)$ with $\nu^{i}(U)<\infty$, we write
\begin{align*}
\widetilde{N}^{i}((0,t]\times U):=N^{i}((0,t]\times U)-t\nu^{i}(U),\quad t\geq 0,
\end{align*}
for the corresponding compensated Poisson random measures on $[0,T]\times \Omega\times \mathbb{R}_0$. See \cite{1989-Ikeda-p555-555} for the details on Poisson random measures.

\vskip 0.3cm

We introduce the following standard space
\[V=\{u\in H^{1}(\mathcal{O})^2: \nabla\cdot u=0, u|_{\partial \mathcal{O}}=0\},\]
with the norm $\|u\|_{V}:=(\int_\mathcal{O}|\nabla u|^2 dx)^{1/2}$ and the inner product $((\cdot,\cdot))$. Denote by $H$ the closure of $V$ in the $L^2$-norm $\|u\|_H:=(\int_\mathcal{O}|u|^2 dx)^{1/2}$. The inner product on $H$ will be denoted by $(\cdot,\cdot)$.

Identifying the Hilbert space $H$ with its dual space $H^*$, via the Riesz representation, we consider the system (\ref{1.a}) in the framework of Gelfand triple:
\begin{align*}
V\subset H\cong H^*\subset V^*.
\end{align*}

\noindent We also denote by $\langle \cdot,\cdot\rangle$ the dual pair between $V^*$ and $V$ from now on.

\vskip 0.3cm

Define the Stokes operator by
\begin{align}\label{A definition}
Au:=-P_H\Delta u,\quad u\in D(A):=H^{2}(\mathcal{O})^2\cap V,
\end{align}
where $P_H:L^2(\mathcal{O})^2\longrightarrow H$ is the usual Helmholtz-Leray projection.
Actually, the map $A$ is an isomorphism between $V$ and $V^*$, and
%$(A,D(A))$ is a self-adjoint operator, moreover,
\begin{align}\label{A self-adjoint}
\langle Au,v\rangle=\langle u,Av\rangle=((u,v)),\quad \forall\,u,v\in V.
\end{align}

\noindent Note that $\|Au\|_{H}$ is a norm on $V\cap H^2(\mathcal{O})^2$ which is equivalent to the Sobolev norm in $H^2(\mathcal{O})^2$(for simplicity denoted by $H^2$ from now on), see Lemma \uppercase\expandafter{\romannumeral3}.3.7 in \cite{1977-Temam-p500-500}.
It is known that there exist an orthonormal basis $\{e_i,i\in\mathbb{N}\}$ in $H$ and corresponding eigenvalues $0<\lambda_i \uparrow<\infty$, that is
\begin{align}\label{A eigenvalue}
Ae_i=\lambda_i e_i,\quad i\in\mathbb{N}.
\end{align}
Since the boundary $\partial \mathcal{O}$ is of class $\mathcal{C}^{3}$, it follows from Chapter \uppercase\expandafter{\romannumeral1}.2.6 in \cite{1977-Temam-p500-500} that
\begin{align}\label{A eigenvector}
e_i\in H^3(\mathcal{O}) .
\end{align}

Set
\begin{align}\label{b definition}
b(u,v,w):=\sum_{i,j}^2\int_{\mathcal{O}}u_i\partial_i v_j w_j\,dx ,\quad u,v,w\in V .
\end{align}
Using integration by parts, it is easy to see that
\begin{align}\label{b eq 0}
b(u,v,w)=-b(u,w,v),\ \  b(u,v,v)=0, \quad u,v,w\in V .
\end{align}

\noindent Throughout the paper, we will denote various generic positive constants by the same letter $C$, although the constants may differ from line to line. We now list some well-known estimates for $b$ which will be used in the sequel(see \cite{1977-Temam-p500-500} for example):
\begin{align}
\label{b estimate 1}&|b(u,v,w)|\leq 2\|u\|_H^{\frac{1}{2}}\|u\|_V^{\frac{1}{2}}\|w\|_H^{\frac{1}{2}}\|w\|_V^{\frac{1}{2}}\|v\|_V, \quad u,v,w\in V , \\
\label{b estimate 2}&|b(u,u,v)|\leq C\|u\|_{H^2}^{\frac{1}{2}}\|u\|_V\|u\|_H^{\frac{1}{2}}\|v\|_{H},\quad u\in V\cap H^2,\quad v\in H.
\end{align}

\noindent For $u,v\in V$, we denote by $B(u,v)$ the element of $V^*$ defined by
\begin{align}\label{B definition}
\langle B(u,v),w\rangle=b(u,v,w),\ \forall\,w\in V.
\end{align}

\noindent Therefore,
\begin{align}\label{B estimate}
\|B(u,v)\|_{V^*}=\sup_{\|w\|_V\leq 1}|b(u,v,w)|\leq 2\|u\|_H^{\frac{1}{2}}\|u\|_V^{\frac{1}{2}}\|v\|_H^{\frac{1}{2}}\|v\|_V^{\frac{1}{2}} ;
\end{align}
hence
\begin{align}
\|B(u,u)\|_{V^*}\leq 2\|u\|_H\|u\|_V.
\end{align}

\noindent We will often use the short notation $B(u):=B(u,u)$. On the other hand, the nonlinear operator $P_H((u\cdot\nabla)v)$ is well defined whenever $u,v$ are such that $(u\cdot\nabla)v$ belongs to $L^2$. One can show that $P_H((u\cdot\nabla)v)$ can be linearly extended to $V\times V\longrightarrow V^*$, and actually coincides with previous $B(u,v)$.

\vskip 0.3cm

It is known that the system (\ref{1.a}) can be reformulated as follows:
\begin{align}\label{Abstract}
\left\{
 \begin{aligned}
 & du(t)=-Au(t)dt-B(u(t),u(t))dt+F(u(t))dt+\sum_{i=1}^m\sigma^i(u(t))dW^i(t), \\
 & u(0)=h .
 \end{aligned}
\right.
\end{align}

%In the following, to simplify notations, we will let $\mu=1$.

%Throughout this paper, we assume that the initial value $\xi$ is measurable with ${\cal F}_0$.
Let $F,\sigma^i, i=1,\cdots,m$ be measurable mappings from $H$ into $H$. We introduce the following condition:

%\hypertarget{12}{text}
%\hyperlink{12}{text}

\noindent \hypertarget{H.1}{{\bf (H.1)}} $F(\cdot), \sigma^i(\cdot) :H\rightarrow H$ are globally Lipschitz maps, i.e.,  there exists a constant $C<\infty$ such that
\begin{align}\label{H1}
\|F(u_1)-F(u_2)\|_H^2 +\sum_{i=1}^m\|\sigma^i (u_1)-\sigma^i (u_2)\|_H^2 \leq C \|u_1-u_2\|_H^2, \quad\forall\, u_1, u_2\in H.
\end{align}

\begin{definition}\label{def u}
A continuous $H$-valued $({\cal F}_t)$-adapted process $u=(u(t))_{t\geq 0}$ is said to be a solution to equation (\ref{Abstract}) if for any $T>0$, $X\in L^2([0,T]\times \Omega, dt\times P, V)$ and for any $t\geq 0$, the following equation holds in $V^*$, $P$-a.s.:
\begin{align}\label{SNSE solution definition}
u(t)=h-\int_0^t Au(s)ds-\int_0^t B(u(s))ds +\int_0^t F(u(s))ds+\sum_{i=1}^m\int_0^t\sigma^i (u(s))dW^i(s) .
\end{align}
\end{definition}
Under the assumption (\hyperlink{H.1}{H.1}) and $h\in H$, it is known that equation (\ref{Abstract})
admits a unique solution (see e.g. \cite{2014-Brzezniak-p283-310}).

\section{Approximations of SNSEs by pure jump type SNSEs}

For $\varepsilon>0$, let $\sigma^{i,\varepsilon}:H\times \mathbb{R}_0\rightarrow H$ be given measurable maps.
Consider the following SNSE driven by pure jump noise:
\begin{align}\label{5.1}
u^{\varepsilon}(t)=& h-\int_0^tAu^{\varepsilon}(s)ds-\int_0^tB(u^{\varepsilon}(s)) ds +\int_0^tF(u^{\varepsilon}(s)) \nonumber\\
&+ \sum_{i=1}^m\int_0^t\int_{\mathbb{R}_0}\sigma^{i,\varepsilon}(u^{\varepsilon}(s-),z)\widetilde{N}^i(dzds).
\end{align}

\noindent We impose the following conditions on $\sigma^{i,\varepsilon}$.

\noindent \hypertarget{H.2}{{\bf (H.2)}} There exists constants $C>0$ and ${\varepsilon}_0>0$ such that
\begin{gather}
\label{H2 H-2} \|F(u)\|_H^2+\sup_{\varepsilon\leq\varepsilon_0}\sum_{i=1}^m \int_{\mathbb{R}_0}\|\sigma^{i,\varepsilon}(u,z)\|_{H}^2\nu^{i}(dz)\leq C(1+\|u\|_H^2), \\
\label{H2 H-4}
\sup_{\varepsilon\leq\varepsilon_0}\sum_{i=1}^m\int_{\mathbb{R}_0}\|\sigma^{i,\varepsilon}(u,z)\|_{H}^4\nu^{i}(dz)\leq C(1+\|u\|_H^4), \\
\label{H2 H Lip}
\|F(u_1)-F(u_2)\|_H^2+\sup_{\varepsilon\leq\varepsilon_0} \sum_{i=1}^m\int_{\mathbb{R}_0}\|\sigma^{i,\varepsilon}(u_1,z)-\sigma^{i,\varepsilon}(u_2,z)\|_{H}^2\nu^{i}(dz)\leq C\|u_1-u_2\|_H^2.
\end{gather}

Denote by $D([0,T], H)$ the space of all c\`a{}dl\`a{}g paths from
$[0,T]$ into $H$ equipped with the Skorohod topology.
\begin{definition}\label{def u varepsilon}
A  $H$-valued $({\cal F}_t)$-adapted process $u^{\varepsilon}=(u^{\varepsilon}(t))_{t\geq 0}$ is said to be a solution to equation (\ref{5.1}) if

(i) for any $T>0$, $u^{\varepsilon}\in D([0,T], H)\cap L^2([0,T]\times \Omega, dt\times P, V)$;

(ii) for every $t\geq 0$, (\ref{5.1}) holds in $V^*$, $P$-a.s..
\end{definition}
Under the assumption (\hyperlink{H.2}{H.2}) and $h\in H$, it is known that for $\varepsilon\leq {\varepsilon}_0$, equation (\ref{5.1}) admits a unique solution (see e.g. \cite{2014-Brzezniak-p283-310}).

%(H2) is a sufficient condition for the existence and uniqueness of the solution of equation (\ref{5.1}), which be found in \cite{2014-Brzezniak-p283-310}.

\vskip 0.3cm

Consider the following conditions.

\noindent \hypertarget{H.3}{{\bf (H.3)}} (i) For each $i\in\{1,\cdots,m\}$, $\forall\,M>0$,
\begin{align}\label{sigma to 0}
\sup_{\|u\|_H\leq M}\sup_{z\in\mathbb{R}_0}\|\sigma^{i,\varepsilon}(u,z)\|_H\xrightarrow{\varepsilon\rightarrow 0} 0.
\end{align}

(ii) For each $i\in\{1,\cdots,m\}$ and each $k,j\in\mathbb{N}$, $u\in H$,
\begin{align}\label{sigma to sigma}
  \int_{\mathbb{R}_0}(\sigma^{i,\varepsilon}(u,z),e_k)(\sigma^{i,\varepsilon}(u,z),e_j)\nu^{i}(dz)\xrightarrow{\varepsilon\rightarrow 0} (\sigma^{i}(u),e_k)(\sigma^{i}(u),e_j).
\end{align}

\noindent \hypertarget{H.4}{{\bf (H.4)}} For each $i\in\{1,\cdots,m\}$ and every $u\in H$,
\begin{align}
\int_{\mathbb{R}_0}\|\sigma^{i,\varepsilon}(u,z)\|_H^2\nu^{i}(dz)\xrightarrow{\varepsilon\rightarrow 0}\|\sigma^{i}(u)\|_H^2.
\end{align}

\begin{remark}
  For our considerations, the jumping measures $\nu^{i}, i=1,\cdots,m$ typically have infinite volume. Therefore, (i),(ii) in (\hyperlink{H.3}{H.3}) and (\hyperlink{H.4}{H.4}) do not imply each other.
\end{remark}

\noindent \hypertarget{H.5}{{\bf (H.5)}} The maps $F,\sigma^{i,\varepsilon}$ take the space $V$ into itself and there exist constants $C>0$ and ${\varepsilon}_0>0$ such that
\begin{align}\label{H3 V}
\|F(u)\|_V^2+\sup_{\varepsilon\leq\varepsilon_0}\sum_{i=1}^m\int_{\mathbb{R}_0}\|\sigma^{i,\varepsilon}(u,z)\|_{V}^2\nu^{i}(dz)\leq C(1+\|u\|_V^2).
\end{align}

%\noindent{\bf (H.3)}$^{\prime}$ There exists a sequence of mappings $\sigma_n^{\varepsilon}(\cdot): H\rightarrow V$ and $F_n(\cdot)H\rightarrow V$ such that for each $n$, $\{\sigma_n^{\varepsilon}\}_{0<\varepsilon <1}$ and $F_n$ satisfies condition {\bf (H.3)}. Moreover,
%
%(i) there is a constant $C$ independent of $n$ such that for every $u\in H$,
%\begin{align}\label{H3 H-2}
%\sup_n\|F_n(u)\|_H^2+\sup_n\|\sigma_n(u)\|_H^2+\sup_{n,\varepsilon}\int_{\mathbb{R}_0}\| \sigma_n^{\varepsilon}(u,z)\|_{H}^2\nu(dz)\leq C(1+\|u\|_H^2).
%\end{align}
%
%
%(ii) there is a constant $C_2$ independent of $n$ such that for every $u_1,u_2\in H$,
%\begin{align}\label{H3 H Lip}
%&\sup_n\|F_n(u_1)-F_n(u_2)\|_H^2+\sup_n\|\sigma_n(u_1)-\sigma_n(u_2)\|_H^2 \nonumber\\ +&\sup_{n,\varepsilon}\int_{\mathbb{R}_0}\| \sigma_n^{\varepsilon}(u_1,z)- \sigma_n^{\varepsilon}(u_2,z)\|_{H}^2\nu(dz)\leq C_2\|u_1-u_2\|_H^2.
%\end{align}
%
%(iii) for every $u\in H$,
%\begin{align}\label{H3 limit}
%\|\sigma_n(u)-\sigma(u)\|_H^2+\sup_{\varepsilon}\int_{\mathbb{R}_0}\| \sigma_n^{\varepsilon}(u,z)-\sigma^{\varepsilon}(u,z)\|_H^2\nu(dz) \xrightarrow{n\rightarrow\infty} 0
%\end{align}

%\begin{remark}
%(H.5) and the (ii) of (H.3) imply that
%\begin{align}
%  \lim_{n\rightarrow\infty}\lim_{\varepsilon\rightarrow 0}\int_{\mathbb{R}_0} \|\sigma_n^{\varepsilon}(x,z)-\sigma^{\varepsilon}(x,z)\|_H^2\nu(dz)=0
%\end{align}
%\end{remark}

\subsection{Preliminary estimates}

We first prepare some preliminary results needed for the proofs of the main results.
In the rest of the paper, for simplicity of the exposition, we let $m=1$ and omit the superscript $i$ of  $\sigma^{i}, \widetilde{N}^{i}, \nu^{i}$.  The case of $m>1$ does not cause  extra difficulties.

\begin{lemma}
Assume (\hyperlink{H.2}{H.2}) and $h\in H$, let $u^{\varepsilon}$ be the solution of equation (\ref{5.1}), then we have
\begin{align}\label{5.0}
\sup_{\varepsilon\leq\varepsilon_0}\Big\{E\sup_{0\leq t\leq T}\|u^{\varepsilon}(t)\|_H^4+E\Big (\int_0^T||u^{\varepsilon}(s)||_V^2ds\Big)^2\Big\}<\infty.
\end{align}
\end{lemma}
\begin{remark}\label{u estimate H}
If we assume (\hyperlink{H.1}{H.1}) and $h\in H$, then (\ref{5.0}) also holds for the solution $u$ of equation (\ref{Abstract}).
\end{remark}

\noindent {\bf Proof}. By It\^{o}'s formula and (\ref{b eq 0}), we have
\begin{align}\label{5.01}
\|u^{\varepsilon}(t)\|_H^2
=&\|h\|_H^2 -2\int_0^t\langle Au^{\varepsilon}(s),u^{\varepsilon}(s)\rangle ds+2 \int_0^t(F(u^{\varepsilon}(s)),u^{\varepsilon}(s))ds\nonumber\\
&+M(t)+ \int_0^t\int_{\mathbb{R}_0}\|\sigma^{\varepsilon}(u^{\varepsilon}(s-),z)\|_H^2 \nu(dz)ds,
\end{align}
where
 \begin{align}\label{martingale M}
 M(t):=\int_0^t\int_{\mathbb{R}_0}\Big( \|\sigma^{\varepsilon}(u^{\varepsilon}(s-),z)\|_H^2+2\big(\sigma^{\varepsilon}(u^{\varepsilon}(s-),z), u^{\varepsilon}(s-)\big)\Big) \widetilde{N}(dzds) .
 \end{align}
Using Burkh\"older's inequality and the assumption (\hyperlink{H.2}{H.2}), we have
 \begin{align}\label{5.02}
&E\sup_{0\leq r\leq t}|M(r)|^2 \nonumber\\
\leq & C E\int_0^t\int_{\mathbb{R}_0}\Big( \|\sigma^{\varepsilon}(u^{\varepsilon}(s-),z)\|_H^2+2\big(\sigma^{\varepsilon}(u^{\varepsilon}(s-),z), u^{\varepsilon}(s-)\big)\Big)^2\nu(dz)ds\nonumber\\
\leq & C E\int_0^t(1+\|u^{\varepsilon}(s)\|_H^4) ds.
\end{align}

\noindent By (\ref{A self-adjoint}), it follows from (\ref{5.01}) that for $t\leq T$,
\begin{align}\label{5.01-1}
\|u^{\varepsilon}(t)\|_H^4 + \Big(\int_0^t\|u^{\varepsilon}(s)\|_V^2ds\Big)^2 \leq C\|h\|_H^4 + C \int_0^t(1+\|u^{\varepsilon}(s)\|_H^4) ds+C M(t)^2.
\end{align}
Take superemum over the interval $[0,t]$ in (\ref{5.01-1}), and use (\ref{5.02}) to get
\begin{align}\label{5.03}
E\sup_{0\leq s\leq t}\|u^{\varepsilon}(s)\|_H^4+E\left (\int_0^t\|X_s^{\varepsilon}\|_V^2ds\right )^2\leq  C \|h\|_H^4 + C E \int_0^t(1+\|u^{\varepsilon}(s)\|_H^4) ds.
\end{align}
Applying Gronwall's inequality completes the proof of the lemma. $\blacksquare$

\begin{lemma}
Assume (\hyperlink{H.2}{H.2}), (\hyperlink{H.5}{H.5}) and $h\in V$. For any constant $M>0$, define
\begin{align}
  \tau_M^{\varepsilon}:=T\wedge\inf\{t\geq 0: \int_0^t\|u^{\varepsilon}(s)\|_V^2ds>M\} \wedge\inf\{t\geq 0: \|u^{\varepsilon}(t)\|_H^2 >M\} ,
\end{align}
where we set $\inf\{\emptyset\}=\infty$. Then we have
\begin{align}\label{V norm estimate}
  \sup_{\varepsilon\leq\varepsilon_0}\Big\{E\sup_{0\leq t\leq \tau_M^{\varepsilon}}\|u^{\varepsilon}(t)\|_V^2+E\Big (\int_0^{\tau_M^{\varepsilon}}\|u^{\varepsilon}(s)\|_{H^2}^2ds\Big)^2\Big\}<\infty.
\end{align}

\end{lemma}

\noindent {\bf Proof}. Through Galerkin approximations, it can be shown that for $\varepsilon\leq\varepsilon_0$, the solution $u^{\varepsilon}\in L^{\infty}([0,T],V)\cap L^2([0,T],H^2)$ with probability one (see e.g. Proposition 2.2 in \cite{2003-Menaldi-p357-381}).
Apply It\^{o}'s formula to $\|u^{\varepsilon}(t)\|_V^2$ to get
\begin{align}\label{u vare V2 eq}
\|u^{\varepsilon}(t)\|_V^2
=&\|h\|_V^2 -2\int_0^t\|Au^{\varepsilon}(s)\|_H^2ds-2\int_0^t\big(B(u^{\varepsilon}(s)),Au^{\varepsilon}(s)\big)ds \nonumber\\
& +2 \int_0^t\big(\big(F(u^{\varepsilon}(s)),u^{\varepsilon}(s)\big)\big)ds +M_1(t)+M_2(t),
\end{align}
where
\begin{align}
 M_1(t):&=2\int_0^t\int_{\mathbb{R}_0}\big(\big(\sigma^{\varepsilon}(u^{\varepsilon}(s-),z), u^{\varepsilon}(s-)\big)\big)\widetilde{N}(dzds), \\
 M_2(t):&=\int_0^t\int_{\mathbb{R}_0}\|\sigma^{\varepsilon}(u^{\varepsilon}(s-),z)\|_V^2 N(dzds).
\end{align}

\noindent Use (\ref{b estimate 2}) and Young's inequality to obtain
\begin{align}\label{BuuAu}
  |(B(u^{\varepsilon}(s)),Au^{\varepsilon}(s))|\leq C\|u^{\varepsilon}\|_{H^2}^{\frac{3}{2}}\|u^{\varepsilon}\|_V\|u^{\varepsilon}\|_H^{\frac{1}{2}} \leq \|u^{\varepsilon}\|_{H^2}^2+C\|u^{\varepsilon}\|_V^4\|u^{\varepsilon}\|_H^2 .
\end{align}

\noindent Therefore, by (\ref{BuuAu}) and (\hyperlink{H.5}{H.5}), we obtain
\begin{align}
&\|u^{\varepsilon}(t)\|_V^2+\int_0^t\|u^{\varepsilon}(s)\|_{H^2}^2ds \nonumber\\
\leq& \|h\|_{V}^2+C\int_0^t\|u^{\varepsilon}(s)\|_V^4\|u^{\varepsilon}(s)\|_H^2ds+ C\int_0^t(1+\|u^{\varepsilon}(s)\|_V^2)ds+M_1(t)+M_2(t).
\end{align}

\noindent Applying Gronwall's inequality yields that
\begin{align}
&\|u^{\varepsilon}(t)\|_V^2+\int_0^t\|u^{\varepsilon}(s)\|_{H^2}^2ds \nonumber\\
\leq& \big(C_T+\|h\|_{V}^2+\sup_{0\leq t\leq \tau^{\varepsilon}_M}|M_1(t)|+M_2(\tau^{\varepsilon}_M)\big) \nonumber\\
&~~~ \times\exp\Big(C_T +C\int_0^t\|u^{\varepsilon}(s)\|_V^2\|u^{\varepsilon}(s)\|_H^2ds\Big),\quad t\in[0,\tau^{\varepsilon}_M].
\end{align}

\noindent Take superemum over the interval $[0,\tau^{\varepsilon}_M]$, remember the definition of $\tau^{\varepsilon}_M$ and take expectations to get
\begin{align}\label{V norm gronwall}
&E\sup_{0\leq t\leq \tau^{\varepsilon}_M}\|u^{\varepsilon}(t)\|_V^2 +E\int_0^{\tau^{\varepsilon}_M}\|u^{\varepsilon}(s)\|_{H^2}^2ds \nonumber\\
\leq& \Big(C_T+\|h\|_{V}^2+E\sup_{0\leq t\leq \tau^{\varepsilon}_M}|M_1(t)|+EM_2(\tau^{\varepsilon}_M)\Big)\exp(C_T+CM^2) .
\end{align}

\noindent By Burkh\"older's inequality, (\hyperlink{H.5}{H.5}) and Young's inequality, we have for $\delta>0$,
\begin{align}\label{EM_1}
E\sup_{0\leq t\leq \tau^{\varepsilon}_M}|M_1(t)|
&\leq 2 E\left[\int_0^{\tau^{\varepsilon}_M}\int_{\mathbb{R}_0} \big(\big(\sigma^{\varepsilon}(u^{\varepsilon}(s-),z), u^{\varepsilon}(s-)\big)\big)^2\nu(dz)ds\right]^{\frac{1}{2}} \nonumber\\
&\leq 2E\left[\int_0^{\tau^{\varepsilon}_M} C\|u^{\varepsilon}(s)\|_V^2 (1+\|u^{\varepsilon}(s)\|_V^2)ds\right]^{\frac{1}{2}} \nonumber\\
&\leq 2(CT+CM)^{\frac{1}{2}} E\sup_{0\leq t\leq \tau^{\varepsilon}_M}\|u^{\varepsilon}(t)\|_V \nonumber\\
&\leq \delta E\sup_{0\leq t\leq \tau^{\varepsilon}_M}\|u^{\varepsilon}(t)\|_V^2+\frac{1}{\delta}(CT+CM) .
\end{align}

\noindent By (\hyperlink{H.5}{H.5}) and (\ref{5.0}), we have
\begin{align}\label{EM_2}
E M_2(\tau^{\varepsilon}_M)
=& E\int_0^{\tau^{\varepsilon}_M}\int_{\mathbb{R}_0}\|\sigma^{\varepsilon}(u^{\varepsilon}(s-),z)\|_V^2 N(dzds) \nonumber\\
=& E\int_0^{\tau^{\varepsilon}_M}\int_{\mathbb{R}_0}\|\sigma^{\varepsilon}(u^{\varepsilon}(s-),z)\|_V^2 \nu(dz)ds \nonumber\\
\leq& E\int_0^{\tau^{\varepsilon}_M}C(1+\|u^{\varepsilon}(s)\|_V^2)ds\leq C<\infty .
\end{align}

\noindent Combining (\ref{V norm gronwall}), (\ref{EM_1}) and (\ref{EM_2}) and choosing sufficiently small $\delta$, we obtain
\begin{align}
E\sup_{0\leq t\leq \tau^{\varepsilon}_M}\|u^{\varepsilon}(t)\|_V^2 +E\int_0^{\tau^{\varepsilon}_M}\|u^{\varepsilon}(s)\|_{H^2}^2ds
\leq C_{T,M}\|h\|_V^2+C_{T,M}
\end{align}
completing the proof of (\ref{V norm estimate}). $\blacksquare$

\vskip 0.3cm

\begin{proposition}\label{un tight on DH}
Assume (\hyperlink{H.2}{H.2}), (\hyperlink{H.5}{H.5}) and $h\in V$. Then the family $\{u^{\varepsilon}, \varepsilon\leq \varepsilon_0\}$ is tight in the space $D([0,T], H)$.
\end{proposition}

\noindent {\bf Proof}. Note that $V$ is compactly embedded into $H$. Thus, by Aldou's tightness criterion (see Theorem 1 in \cite{1978-Aldous-p335-340}), it suffices to show that:

\noindent (i) for any $0<\eta<1$, there exists $L_{\eta}>0$ such that
\begin{align}\label{i verify compact}
  \sup_{\varepsilon\leq\varepsilon_0}P\left(\sup_{0\leq t\leq T}\|u^{\varepsilon}(t)\|_V>L_{\eta}\right)<\eta ;
\end{align}

\noindent (ii) for any stopping time $0\leq\zeta^{\varepsilon}\leq T$ with respect to the natural filtration generated by $\{u^{\varepsilon}(s), s\leq t\}$, and any $\eta>0$,
%and any constant sequence $\delta^{\varepsilon}$ convergence to $0$,
%\begin{align}\label{ii verfity conti}
%\sup_{\varepsilon}E\sup_{\tau^{\varepsilon}\leq r\leq (\tau^{\varepsilon}+\theta)\wedge T}\|u^{\varepsilon}(r)-u^{\varepsilon}(\tau^{\varepsilon})\|_H^2\leq C\theta+C\theta^{\frac{1}{2}},\quad \forall\,\theta>0.
%\end{align}
\begin{align}\label{ii verfity conti}
\lim_{\delta\rightarrow0}\sup_{\varepsilon\leq\varepsilon_0}P(\|u^{\varepsilon}(\zeta^{\varepsilon}+\delta) -u^{\varepsilon}(\zeta^{\varepsilon})\|_H>\eta)=0 ,
\end{align}
%\begin{align}\label{ii verfity conti}
%u^{\varepsilon}(\zeta^{\varepsilon}+\delta^{\varepsilon})-u^{\varepsilon}(\zeta^{\varepsilon}) \xrightarrow{\varepsilon\rightarrow 0} 0\quad \text{in probability}.
%\end{align}
where we set $\zeta^{\varepsilon}+\delta :=T\wedge(\zeta^{\varepsilon}+\delta)$.
%Since (i) implies that, for each $t\in[0,T]$, $\{u^{\varepsilon}(t)\}$ is tight on $D([0,T],H)$;

\vskip 0.3cm

Note that (\ref{5.0}) implies
\begin{align}\label{tau less T}
 & \sup_{\varepsilon\leq\varepsilon_0}P\left(\tau_M^{\varepsilon}<T\right)\nonumber\\
\leq & \sup_{\varepsilon\leq\varepsilon_0}P\left(\int_0^T\|u^{\varepsilon}(s)\|_V^2ds>M\right) +\sup_{\varepsilon\leq\varepsilon_0}P\left(\sup_{0\leq s\leq T}\|u^{\varepsilon}(s)\|_H^2>M\right)\nonumber\\
\leq & \frac{1}{M}\sup_{\varepsilon\leq\varepsilon_0}E\int_0^T\|u^{\varepsilon}(s)\|_V^2ds +\frac{1}{M}\sup_{\varepsilon\leq\varepsilon_0} E\sup_{0\leq t\leq T}\|u^{\varepsilon}(s)\|_H^2\nonumber\\
\leq & \frac{C}{M} .
\end{align}

\noindent For any $L>0$, by (\ref{tau less T}) and (\ref{V norm estimate}), we have
\begin{align}\label{52}
&\sup_{\varepsilon\leq\varepsilon_0}P\left(\sup_{0\leq t\leq T}\|u^{\varepsilon}(t)\|_V>L\right)\nonumber\\
\leq& \sup_{\varepsilon\leq\varepsilon_0}P\left(\sup_{0\leq t\leq T}\|u^{\varepsilon}(t)\|_V>L, \tau_M^{\varepsilon}=T\right)+\sup_{\varepsilon\leq\varepsilon_0}P\left(\tau_M^{\varepsilon}<T\right)\nonumber\\
\leq& \sup_{\varepsilon\leq\varepsilon_0}P\left(\sup_{0\leq t\leq \tau_M^{\varepsilon}}\|u^{\varepsilon}(t)\|_V>L\right) + \frac{C}{M}\nonumber\\
\leq & \frac{1}{L^2}\sup_{\varepsilon\leq\varepsilon_0}E\sup_{0\leq t\leq \tau_M^{\varepsilon}}\|u^{\varepsilon}(t)\|_V^2 +\frac{C}{M}\nonumber\\
\leq & \frac{C_M}{L^2}+\frac{C}{M} .
\end{align}

\noindent Given any $\eta>0$, we can first take sufficiently large constant $M$, and then choose the constant $L$ so that the right hand side of (\ref{52}) will be smaller than $\eta$. Hence (i) is satisfied.

Now, we come to verify (ii). For any $\eta>0$,
\begin{align}
 &\sup_{\varepsilon\leq\varepsilon_0} P(\|u^{\varepsilon}(\zeta^{\varepsilon}+\delta)-u^{\varepsilon}(\zeta^{\varepsilon})\|_H>\eta) \nonumber\\
\leq & \sup_{\varepsilon\leq\varepsilon_0} P\Big(\Big\|\int_{\zeta^{\varepsilon}}^{\zeta^{\varepsilon}+\delta}Au^{\varepsilon}(s)ds\Big\|_H>\frac{\eta}{4}\Big) \nonumber\\
& + \sup_{\varepsilon\leq\varepsilon_0} P\Big(\Big\|\int_{\zeta^{\varepsilon}}^{\zeta^{\varepsilon}+\delta}B(u^{\varepsilon}(s))ds\Big\|_H>\frac{\eta}{4}\Big) \nonumber\\
& + \sup_{\varepsilon\leq\varepsilon_0} P\Big(\Big\|\int_{\zeta^{\varepsilon}}^{\zeta^{\varepsilon}+\delta}F(u^{\varepsilon}(s))ds\Big\|_H>\frac{\eta}{4}\Big) \nonumber\\
& +\sup_{\varepsilon\leq\varepsilon_0} P\Big(\Big\|\int_{\zeta^{\varepsilon}}^{\zeta^{\varepsilon}+\delta}\int_{\mathbb{R}_0} \sigma^{\varepsilon}(u^{\varepsilon}(s-),z)\widetilde{N}(dzds)\Big\|_H>\frac{\eta}{4}\Big) \nonumber\\
& :=I_1+I_2+I_3+I_4 .
\end{align}

\noindent By H\"older's inequality and Chebyshev's inequality, it follows from (\ref{V norm estimate}) and (\ref{tau less T}) that for $M>0$,
\begin{align}\label{I1 delta to 0}
 I_1\leq & \sup_{\varepsilon\leq\varepsilon_0} P\Big(\delta\int_{\zeta^{\varepsilon}}^{\zeta^{\varepsilon}+\delta} \|Au^{\varepsilon}(s)\|_H^2ds>\frac{\eta^2}{16}\Big) \nonumber\\
\leq & \sup_{\varepsilon\leq\varepsilon_0}P\Big(\delta\int_{\zeta^{\varepsilon}}^{\zeta^{\varepsilon}+\delta} \|Au^{\varepsilon}(s)\|_H^2ds>\frac{\eta^2}{16},\tau^{\varepsilon}_M =T\Big) + \sup_{\varepsilon\leq\varepsilon_0}P(\tau^{\varepsilon}_M<T) \nonumber\\
\leq &\sup_{\varepsilon\leq\varepsilon_0} P\Big(\delta\int_0^{\tau^{\varepsilon}_M} \|Au^{\varepsilon}(s)\|_H^2ds>\frac{\eta^2}{16}\Big) + \frac{C}{M} \nonumber\\
\leq & \frac{16}{\eta^2}\delta\sup_{\varepsilon\leq\varepsilon_0}E\int_0^{\tau^{\varepsilon}_M} \|Au^{\varepsilon}(s)\|_H^2ds + \frac{C}{M} \nonumber\\
\leq & \frac{C_M}{\eta^2} \delta + \frac{C}{M} .
\end{align}

\noindent By (\ref{b estimate 2}), we have
$\|B(u^{\varepsilon}(s))\|_H\leq C\|u^{\varepsilon}(s)\|_{H^2}^{\frac{1}{2}}\|u^{\varepsilon}(s)\|_V\|u^{\varepsilon}(s)\|_H^{\frac{1}{2}}$.
Using (\ref{V norm estimate}) and (\ref{tau less T}), we have
\begin{align}\label{I2 delta to 0}
I_2\leq & \sup_{\varepsilon\leq\varepsilon_0} P\Big(\int_{\zeta^{\varepsilon}}^{\zeta^{\varepsilon}+\delta}\|B(u^{\varepsilon}(s))\|_H ds>\frac{\eta}{4}\Big) \nonumber\\
\leq & \sup_{\varepsilon\leq\varepsilon_0} P\Big(\int_{\zeta^{\varepsilon}}^{\zeta^{\varepsilon}+\delta} \|u^{\varepsilon}(s)\|_{H^2}^{\frac{1}{2}}\|u^{\varepsilon}(s)\|_V\|u^{\varepsilon}(s)\|_H^{\frac{1}{2}} ds>\frac{\eta}{4C}\Big) \nonumber\\
\leq & \sup_{\varepsilon\leq\varepsilon_0} P\Big(\int_{\zeta^{\varepsilon}}^{\zeta^{\varepsilon}+\delta} \|u^{\varepsilon}(s)\|_{H^2}^{\frac{1}{2}}\|u^{\varepsilon}(s)\|_V\|u^{\varepsilon}(s)\|_H^{\frac{1}{2}} ds>\frac{\eta}{4C}, \tau^{\varepsilon}_M =T\Big) \nonumber\\
& +  \sup_{\varepsilon\leq\varepsilon_0}P(\tau^{\varepsilon}_M<T) \nonumber\\
\leq & \sup_{\varepsilon\leq\varepsilon_0} P\Big(\int_{\zeta^{\varepsilon}}^{(\zeta^{\varepsilon}+\delta)\wedge\tau^{\varepsilon}_M} \|u^{\varepsilon}(s)\|_{H^2}^{\frac{1}{2}}\|u^{\varepsilon}(s)\|_V\|u^{\varepsilon}(s)\|_H^{\frac{1}{2}} ds>\frac{\eta}{4C}\Big) + \frac{C}{M} \nonumber\\
\leq & \frac{4C}{\eta} \sup_{\varepsilon\leq\varepsilon_0} \bigg[ \Big(E\int_{\zeta^{\varepsilon}}^{(\zeta^{\varepsilon}+\delta)\wedge\tau^{\varepsilon}_M} \|u^{\varepsilon}(s)\|_{H}^2ds\Big)^\frac{1}{4}
\Big(E\int_{\zeta^{\varepsilon}}^{(\zeta^{\varepsilon}+\delta)\wedge\tau^{\varepsilon}_M} \|u^{\varepsilon}(s)\|_{H^2}^2ds\Big)^\frac{1}{4} \nonumber\\
& \times\Big(E\int_{\zeta^{\varepsilon}}^{(\zeta^{\varepsilon}+\delta)\wedge\tau^{\varepsilon}_M} \|u^{\varepsilon}(s)\|_V^2ds\Big)^\frac{1}{2} \bigg] + \frac{C}{M} \nonumber\\
\leq & \frac{C_M}{\eta}\delta^{\frac{3}{4}} \sup_{\varepsilon\leq\varepsilon_0}\Big(E\sup_{0\leq s\leq T}\|u^{\varepsilon}(s)\|_H^2\Big)^{\frac{1}{4}}\times
\sup_{\varepsilon\leq\varepsilon_0} \Big(E\int_0^{\tau^{\varepsilon}_M}\|u^{\varepsilon}(s)\|_{H^2}^2ds\Big)^{\frac{1}{4}} \nonumber\\
& \times \sup_{\varepsilon\leq\varepsilon_0} \Big(E\sup_{0\leq s\leq \tau^{\varepsilon}_M}\|u^{\varepsilon}(s)\|_V^2\Big)^{\frac{1}{2}}+\frac{C}{M} \nonumber\\
\leq & \frac{C_M}{\eta}\delta^{\frac{3}{4}} +\frac{C}{M} .
\end{align}

\noindent On the other hand, by (\hyperlink{H.2}{H.2}) and (\ref{5.0}) we have
\begin{align}\label{I3 delta to 0}
 I_3 \leq & \frac{4}{\eta}\sup_{\varepsilon\leq\varepsilon_0} E\int_{\zeta^{\varepsilon}}^{\zeta^{\varepsilon}+\delta}\|F(u^{\varepsilon}(s))\|_H ds \nonumber\\
 \leq & \frac{4}{\eta}\sup_{\varepsilon\leq\varepsilon_0} E\int_{\zeta^{\varepsilon}}^{\zeta^{\varepsilon}+\delta} C(1+\|u^{\varepsilon}(s)\|_H)ds \nonumber\\
 \leq & \frac{C}{\eta}\delta \Big(1+\sup_{\varepsilon\leq\varepsilon_0} E\sup_{0\leq s\leq T}\|u^{\varepsilon}(s)\|_H\Big) \nonumber\\
 \leq & \frac{C}{\eta}\delta .
\end{align}

\noindent Similarly,
\begin{align}\label{I4 delta to 0}
I_4 \leq & \frac{16}{\eta^2}\sup_{\varepsilon\leq\varepsilon_0} E\Big\|\int_{\zeta^{\varepsilon}}^{\zeta^{\varepsilon}+\delta}\int_{\mathbb{R}_0} \sigma^{\varepsilon}(u^{\varepsilon}(s-),z)\widetilde{N}(dzds)\Big\|_H^2 \nonumber\\
\leq & \frac{16}{\eta^2}\sup_{\varepsilon\leq\varepsilon_0} E\int_{\zeta^{\varepsilon}}^{\zeta^{\varepsilon}+\delta}\int_{\mathbb{R}_0} \|\sigma^{\varepsilon}(u^{\varepsilon}(s-),z)\|_H^2\nu(dz)ds \nonumber\\
\leq & \frac{16}{\eta^2}\sup_{\varepsilon\leq\varepsilon_0} E\int_{\zeta^{\varepsilon}}^{\zeta^{\varepsilon}+\delta} C(1+\|u^{\varepsilon}(s)\|_H^2)ds \nonumber\\
%\leq & \frac{C}{\eta^2}\delta \Big(1+\sup_{\varepsilon} E\sup_{0\leq s\leq T}\|u^{\varepsilon}(s)\|_H^2\Big) \nonumber\\
\leq & \frac{C}{\eta^2}\delta .
\end{align}

\noindent Combine (\ref{I1 delta to 0})---(\ref{I4 delta to 0}) together,
first let $\delta\rightarrow 0$, then let $M\rightarrow\infty$ to obtain (\ref{ii verfity conti}).
Thus (ii) is verified, which completes the proof. $\blacksquare$

\subsection{The weak convergence}

Denote by $\mu_{\varepsilon}$, $\mu$ respectively the laws of $u^{\varepsilon}$ and $u$ on the spaces $D([0,T], H)$ and $C([0,T],H)$. We will establish the weak convergence by two stages. We first obtain the weak convergence in Theorem 3.7  under stronger conditions, and then we remove the extra assumptions and get the general convergence result in Theorem 3.8.
\begin{theorem}\label{V growth theorem}
Assume (\hyperlink{H.1}{H.1}), (\hyperlink{H.2}{H.2}), (\hyperlink{H.3}{H.3}), (\hyperlink{H.5}{H.5}) and $h\in V$. Then,
for any $T>0$, $\mu_{\varepsilon}$ converges weakly to $\mu$, as $\varepsilon\rightarrow 0$, on the space $D([0,T], H)$ equipped with the Skorohod topology.
\end{theorem}

\noindent {\bf Proof}. By Proposition \ref{un tight on DH}, the family $\{ \mu_{\varepsilon}, \varepsilon\leq\varepsilon_0\}$ is tight in $D([0,T], H)$. Let $ \mu_0$ be the weak limit of any convergent subsequence $\{\mu_{\varepsilon_n}\}$. We will show that $\mu_0=\mu$.
The rest of the proof is divided into three steps. In step 1, we show that $\mu_0$ is supported on the space $C([0,T],H)$. In step 2, we prove that $\mu_0$ is a solution of a martingale problem. In step 3, we show that $\mu_0$ is the law of a weak solution of SNSE (\ref{Abstract}), hence complete the proof.

\vskip 0.3cm

Step 1. For any $\eta>0, M>0$, we have
\begin{align}
&P\left(\sup_{0< t\leq T}\|u^{\varepsilon}(t)-u^{\varepsilon}(t-)\|_H>\eta\right) \nonumber\\
\leq & P\left(\sup_{0< t\leq T}\sup_{z\in\mathbb{R}_0}\|\sigma^{\varepsilon}(u^{\varepsilon}(t-),z)\|_H>\eta\right)\nonumber\\
\leq & P\left(\sup_{0\leq t\leq T}\sup_{z\in\mathbb{R}_0}\|\sigma^{\varepsilon}(u^{\varepsilon}(t),z)\|_H>\eta, \sup_{0\leq t\leq T}\|u^{\varepsilon}(t)\|\leq M\right) \nonumber\\
&+ P\left(\sup_{0\leq t\leq T}\|u^{\varepsilon}(t)\|> M\right)\nonumber\\
\leq & P\Big(\sup_{\|x\|_H\leq M}\sup_{z\in\mathbb{R}_0}\|\sigma^{\varepsilon}(x,z)\|_H>\eta\Big) +\frac{1}{M^2}\sup_{\varepsilon\leq\varepsilon_0}E\sup_{0\leq t\leq T}\|u^{\varepsilon}(t)\|_H^2 .
\end{align}

\noindent By (\ref{5.0}) and (\ref{sigma to 0}), we first let $\varepsilon\rightarrow 0$ and then $M\rightarrow \infty$ to see that
\begin{align}
\sup_{0< t\leq T}\|u^{\varepsilon}(t)-u^{\varepsilon}(t-)\|_H\xrightarrow{\varepsilon\rightarrow 0} 0 \quad\text{in probability}.
\end{align}

\noindent Therefore, it follows from Theorem 13.4 in \cite{1999-Billingsley-p277-277} that $\mu_0$ is supported on the space $C([0,T],H)$.  As a consequence, the finite dimensional distributions of $\mu_{\varepsilon_n}$ converge to that of $\mu_0$.

\vskip 0.3cm

Step 2. For $k,j\in\mathbb{N}$, let $f(x)=(x,e_k)(x,e_j)$, $x\in H$.
The gradient of $f$ (denoted by $\nabla f$) and the operator (denoted by $f^{\prime\prime}$) associated with the second derivatives of $f$ are respectively given by
\begin{gather}
\nabla f(x)=(x,e_j)e_k+(x,e_k)e_j ,\\
f^{\prime\prime}(x)=e_j\otimes e_k+e_k\otimes e_j.
\end{gather}

\noindent Set
\begin{align}
\label{definition Lvare}
L^{\varepsilon}f(x):=&-(A\nabla f(x),x)-\langle B(x),\nabla f(x)\rangle +(F(x),\nabla f(x)) \nonumber\\ &+\int_{\mathbb{R}_0}\big[f(x+\sigma^{\varepsilon}(x,z))-f(x)-(\nabla f(x),\sigma^{\varepsilon}(x,z))\big]\nu(dz), \\
\label{definition L}
Lf(x):=&-(A\nabla f(x),x)-\langle B(x),\nabla f(x)\rangle +(F(x),\nabla f(x)) +\frac{1}{2}(f^{\prime\prime}(x)\sigma(x),\sigma(x)) .
\end{align}

\noindent By It\^{o}'s formula,
\begin{align}
&f(u^{\varepsilon}(t))-f(h)-\int_0^t L^{\varepsilon}f(u^{\varepsilon}(s))ds\nonumber\\
=&\int_0^t\int_{\mathbb{R}_0}\Big[f\big(u^{\varepsilon}(s-)+\sigma^{\varepsilon}(u^{\varepsilon}(s-),z)\big) -f\big(u^{\varepsilon}(s-)\big)\Big]\widetilde{N}(dzds)
\end{align}
is a martingale.
Denote by $X_t(\omega):=\omega(t)$, $\omega\in D([0,T],H)$ the coordinate process on $D([0,T],H)$. By the above martingale property, for any $0\leq s_0< s_1<...<s_n\leq s<t$ and $f_0, f_1, ...f_n\in C_b(H)$(the collection of bounded  continuous functions on $H$), it holds that
\begin{align}\label{5.23}
E^{\mu_{\varepsilon}}\left[\left(f(X_t)-f(X_s)-\int_s^tL^{\varepsilon}f(X_r)dr\right) f_0(X_{s_0})...f_n(X_{s_n})\right] = 0.
\end{align}

\noindent Let
\begin{align}
G_{\varepsilon}(x):=\Big|\int_{\mathbb{R}_0} (\sigma^{\varepsilon}(x,z),e_k)(\sigma^{\varepsilon}(x,z),e_j)\nu(dz)-(\sigma(x),e_k)(\sigma(x),e_j)\Big| ,
\end{align}

\noindent $x\in H$. By (\ref{definition Lvare}) and (\ref{definition L}), we have
\begin{align}\label{Lvarepsilon-L}
|L^{\varepsilon}f(X_r)-Lf(X_r)|=G_{\varepsilon}(X_r) .
\end{align}

\noindent We claim that
\begin{align}\label{5.26}
\lim_{n\rightarrow \infty}E^{\mu_{\varepsilon_n}}\Big[\int_s^t|L^{\varepsilon_n}f(X_r)-Lf(X_r)|dr\Big] = 0 . %\lim_{n\rightarrow \infty}\int_s^t EG_{\varepsilon_n}(u^{\varepsilon_n}(r))dr=0  .
\end{align}

\noindent Note that
\begin{gather}
E^{\mu_{\varepsilon_n}}\Big[\int_s^t|L^{\varepsilon_n}f(X_r)-Lf(X_r)|dr\Big] =\int_s^t E G_{\varepsilon_n}(u^{\varepsilon_n}(r))dr,\\
\label{G growth}  \sup_{\varepsilon\leq\varepsilon_0}G_{\varepsilon}(x)\leq C(1+\|x\|_H^2).
\end{gather}

\noindent By the dominated convergence theorem and (\ref{5.0}), to prove (\ref{5.26}),
%\begin{align}\label{G growth}
%  \sup_{\varepsilon}G_{\varepsilon}(x)\leq C(1+\|x\|_H^2),
%\end{align}
it suffices to prove that for every $r\in[0,T]$,
\begin{align}\label{EG varen to 0}
\lim_{n\rightarrow\infty}EG_{\varepsilon_n}(u^{\varepsilon_n}(r))= 0 .
\end{align}

\noindent Since the finite dimensional distributions of $\mu_{\varepsilon_n}$ converge weakly to that of $\mu_0$, by the Skorohod's representation theorem, in order not to introduce more notations, we can assume that $u^{\varepsilon_n}(r)$ converges almost surely to a $H$-valued random variable $u^0$. In view of (\ref{5.0}), $\{\|u^{\varepsilon_n}(r)\|_H^2\}_{n\geq 1}$ is uniformly integrable, and therefore we can further deduce that $u^0\in L^2(\Omega,H)$ and
\begin{align}\label{E uvare-u0 to 0}
\lim_{n\rightarrow\infty}E\|u^{\varepsilon_n}(r)-u^0\|_H^2= 0 .
\end{align}

\noindent By the dominated convergence theorem, it follows from (\ref{sigma to sigma}) and (\ref{G growth}) that
\begin{align}\label{EG u0 to 0}
\lim_{n\rightarrow\infty}EG_{\varepsilon_n}(u^{0})=0.
\end{align}
Hence to prove (\ref{EG varen to 0}), it suffices to prove
%\noindent Now, (\ref{EG varen to 0}) is followed by the the following equation:
%\begin{align}\label{EG u0 to 0}
% \lim_{n\rightarrow\infty}EG_{\varepsilon_n}(u^{\varepsilon_n}(r))= \lim_{n\rightarrow\infty}EG_{\varepsilon_n}(u^{0})=0.
%\end{align}
%
%\noindent The second equal sign in (\ref{EG u0 to 0}) is established by (\ref{sigma to sigma}), (\ref{G growth}) and the dominated convergence theorem. To prove the first equal sign, it suffice to show that
\begin{align}\label{G-G to 0}
\lim_{n\rightarrow\infty}E|G_{\varepsilon_n}(u^{\varepsilon_n}(r))-G_{\varepsilon_n}(u^0)|=0 .
\end{align}
We have
\begin{align}
&E|G_{\varepsilon_n}(u^{\varepsilon_n}(r))-G_{\varepsilon_n}(u^0)| \nonumber\\
\leq & E\Big|\int_{\mathbb{R}_0} \big(\sigma^{\varepsilon_n}(u^{\varepsilon_n}(r),z),e_k\big) \big(\sigma^{\varepsilon_n}(u^{\varepsilon_n}(r),z),e_j\big)\nu(dz) \nonumber\\
&~~~-\int_{\mathbb{R}_0}\big(\sigma^{\varepsilon_n}(u^0,z),e_k\big) \big(\sigma^{\varepsilon_n}(u^0,z),e_j\big)\nu(dz)\Big| \nonumber\\
& +E|(\sigma(u^{\varepsilon_n}(r)),e_k)(\sigma(u^{\varepsilon_n}(r)),e_j) -(\sigma(u^0),e_k)(\sigma(u^0),e_j)|\nonumber\\
:=& I_1+I_2 .
\end{align}

\noindent In view of (\ref{H2 H-2}) and (\ref{H2 H Lip}), we have
\begin{align}
  I_1\leq & E\int_{\mathbb{R}_0}\Big| \big(\sigma^{\varepsilon_n}(u^{\varepsilon_n}(r),z),e_k\big)\big(\sigma^{\varepsilon_n}(u^{\varepsilon_n}(r),z) -\sigma^{\varepsilon_n}(u^0,z),e_j\big)\Big|\nu(dz) \nonumber\\
&+E\int_{\mathbb{R}_0}\Big|\big(\sigma^{\varepsilon_n}(u^{\varepsilon_n}(r),z)-\sigma^{\varepsilon_n}(u^0,z),e_k\big) \big(\sigma^{\varepsilon_n}(u^0,z),e_j\big)\Big|\nu(dz) \nonumber\\
\leq & \Big[E\int_{\mathbb{R}_0} \big\|\sigma^{\varepsilon_n}(u^{\varepsilon_n}(r),z)\big\|_H^2\nu(dz)\Big]^{\frac{1}{2}} \Big[E\int_{\mathbb{R}_0}\big\|\sigma^{\varepsilon_n}(u^{\varepsilon_n}(r),z) -\sigma^{\varepsilon_n}(u^0,z)\big\|_H^2\nu(dz)\Big]^{\frac{1}{2}} \nonumber\\
&+\Big[E\int_{\mathbb{R}_0} \big\|\sigma^{\varepsilon_n}(u^0,z)\big\|_H^2\nu(dz)\Big]^{\frac{1}{2}} \Big[E\int_{\mathbb{R}_0}\big\|\sigma^{\varepsilon_n}(u^{\varepsilon_n}(r),z) -\sigma^{\varepsilon_n}(u^0,z)\big\|_H^2\nu(dz)\Big]^{\frac{1}{2}}\nonumber\\
\leq & C\Big[(1+E\|u^0\|_H^2)^{\frac{1}{2}}+\sup_{\varepsilon_n}(1+E\|u^{\varepsilon_n}(r)\|_H^2)^{\frac{1}{2}}\Big] \big(E\|u^{\varepsilon_n}(r)-u^0\|_H^2\big)^{\frac{1}{2}} .
\end{align}

\noindent This yields that $I_1\rightarrow 0$ taking into account (\ref{5.0}) and (\ref{E uvare-u0 to 0}). A similar argument leads to  $I_2\rightarrow 0$. Therefore, (\ref{G-G to 0}) holds. Hence the claim (\ref{5.26}) is proved.

%\begin{align}\label{5.26}
%&\lim_{n\rightarrow \infty}E^{\mu_{\varepsilon_n}}\int_s^t|L^{\varepsilon_n}f(\Pi_r)-Lf(\Pi_r)|dr \nonumber\\
%=& \int_s^t E\Big|\int_{\mathbb{R}_0} (\sigma^{\varepsilon_n}(u^{\varepsilon_n}(r),z),e_k)(\sigma^{\varepsilon_n}(u^{\varepsilon_n}(r),z),e_j)\nu(dz)\nonumber\\ &~~~~~~~~~~-(\sigma(u^{\varepsilon_n}(r)),e_k)(\sigma(u^{\varepsilon_n}(r)),e_j)\Big|dr \rightarrow 0 .
%\end{align}

%\begin{align}\label{5.26}
%&\lim_{\varepsilon\rightarrow 0}E^{\mu_{\varepsilon}}\int_s^t|L^{\varepsilon}f(\Pi_r)-Lf(\Pi_r)|dr \nonumber\\
%=&\int_s^t E^{\mu_{\varepsilon}}\Big|\int_{\mathbb{R}_0} (\sigma^{\varepsilon}(\Pi_r,z),e_k)(\sigma^{\varepsilon}(\Pi_r,z),e_j)\nu(dz)\nonumber\\
%&~~~~~~~~~~-(\sigma(\Pi_r),e_k)(\sigma(\Pi_r),e_j)\Big|dr \nonumber\\
%=& \int_s^t E^{\mu_{\varepsilon}\Pi_r^{-1}}\Big|\int_{\mathbb{R}_0} (\sigma^{\varepsilon}(x,z),e_k)(\sigma^{\varepsilon}(x,z),e_j)\nu(dz) \nonumber\\ &~~~~~~~~~~~~~~~-(\sigma(x),e_k)(\sigma(x),e_j)\Big|dr
%\end{align}

\vskip 0.3cm
Next we prove that
\begin{align}\label{Mkj martingale}
 M_{k,j}(t):=f(X_t)-f(h)-\int_0^tLf(X_r)dr
\end{align}
is a martingale under $\mu_0$. This is equivalent to proving that
\begin{align}\label{5.27}
 E^{\mu_{0}}\Big[\Big(f(X_t)-f(X_s)-\int_s^t Lf(X_r)dr\Big)f_0(X_{s_0})...f_n(X_{s_n})\Big] =0 .
\end{align}

\noindent Since the finite dimensional distributions of $\mu_{\varepsilon_n}$ converge to that of $\mu_0$, noticing that $\|f(x)\|_H\leq\|x\|_H^2$ and (\ref{5.0}), it follows from Theorem 1.6.8 in \cite{2010-Durrett-p428-428} that
\begin{align}\label{E vareF to EF}
E^{\mu_{0}}\Big[f(X_t)f_0(X_{s_0})\cdots f_n(X_{s_n})\Big] =\lim_{n\rightarrow \infty}E^{\mu_{\varepsilon_n}}\Big[f(X_t)f_0(X_{s_0})...f_n(X_{s_n})\Big] .
\end{align}

%To show
%\begin{align}
%  &E^{\mu_{0}}\Big[\Big(\int_s^t Lf(\Pi_r)dr\Big)f(\Pi_{s_0})...f(\Pi_{s_n})\Big] \nonumber\\
%=&\lim_{n\rightarrow \infty}E^{\mu_{\varepsilon_n}}\Big[\Big(\int_s^t Lf(\Pi_r)dr\Big)f(\Pi_{s_0})...f(\Pi_{s_n})\Big] \nonumber\\
%\end{align}
%Since
%\begin{align}
%  &E^{\mu_{0}}\Big[\Big(\int_s^t Lf(\Pi_r)dr\Big)f(\Pi_{s_0})...f(\Pi_{s_n})\Big] \nonumber\\
%=&\int_s^t E^{\mu_{\varepsilon_n}}\Big[ \Big(Lf(\Pi_r)\Big)f(\Pi_{s_0})...f(\Pi_{s_n})\Big]dr \nonumber\\
%\end{align}

\noindent In view of (\ref{A eigenvector}), we have
\begin{align}
|\langle B(x,x),e_k\rangle|=|\langle B(x,e_k),x\rangle|\leq C\|x\|_H^2\|\nabla e_k\|_{L^{\infty}}\leq C\|e_k\|_{H^3}\|x\|_H^2 .
\end{align}
Thus, $Lf(x)$ is a continuous function on $H$ and
\begin{align}
|Lf(x)|\leq C(1+\|x\|_H^3) .
\end{align}
Therefore, for the same reason as (\ref{E vareF to EF}), we have for every $r\in[s,t]$,
\begin{align}\label{varepsilon L to 0 L}
  E^{\mu_{0}}\big[\big(Lf(X_r)\big)f_0(X_{s_0})...f_n(X_{s_n})\big]
=\lim_{n\rightarrow \infty}E^{\mu_{\varepsilon_n}}\big[\big(Lf(X_r)\big)f_0(X_{s_0})...f_n(X_{s_n})\big].
\end{align}

\noindent By the Fubini theorem and the dominate convergence theorem, we obtain
\begin{align}\label{Emun of L to Emu0}
  &E^{\mu_{0}}\Big[\Big(\int_s^t Lf(X_r)dr\Big)f_0(X_{s_0})...f_n(X_{s_n})\Big] \nonumber\\
=&\lim_{n\rightarrow \infty}E^{\mu_{\varepsilon_n}}\Big[\Big(\int_s^t Lf(X_r)dr\Big)f_0(X_{s_0})...f_n(X_{s_n})\Big] .
\end{align}

\noindent Using (\ref{E vareF to EF}), (\ref{Emun of L to Emu0}), (\ref{5.26}) and (\ref{5.23}), we have
\begin{align}%\label{5.27}
& E^{\mu_{0}}\Big[\Big(f(X_t)-f(X_s)-\int_s^t Lf(X_r)dr\Big)f_0(X_{s_0})...f_n(X_{s_n})\Big] \nonumber\\
=&\lim_{n\rightarrow \infty}E^{\mu_{\varepsilon_n}}\Big[\Big(f(X_t)-f(X_s)-\int_s^t Lf(X_r)dr\Big)f_0(X_{s_0})...f_n(X_{s_n})\Big] \nonumber\\
=&\lim_{n\rightarrow \infty}E^{\mu_{\varepsilon_n}}\Big[\Big(f(X_t)-f(X_s)-\int_s^t L^{\varepsilon_n}f(X_r)dr\Big)f_0(X_{s_0})...f_n(X_{s_n})\Big]\nonumber\\
=& 0 .
\end{align}

\noindent Hence $M_{k,j}(t)$ in (\ref{Mkj martingale}) is a martingale under $\mu_0$.

%
%
%\noindent This establishes the first equal sign in (\ref{5.27}), hence (\ref{5.27}) holds.
%
%Since $0\leq s_0< s_1<...<s_n\leq s<t$ are arbitrary and $f(x)=(x,e_k)(x,e_j)$, (\ref{5.27}) implies that
%\begin{align}\label{Mkj martingale}
%& M_{k,j}(t):=f(\Pi_t)-f(\xi)-\int_0^tLf(\Pi_r)dr \nonumber\\
%=& (\Pi_t,e_k)(\Pi_t,e_j)-(\xi,e_k)(\xi,e_j)\nonumber\\
%&+\int_0^t\big[(Ae_k,\Pi_s)(\Pi_s,e_j)+(Ae_j,\Pi_s)(\Pi_s,e_k)\big]ds\nonumber\\
%&+\int_0^t\langle B(\Pi_s),(\Pi_s,e_j)e_k+(\Pi_s,e_k)e_j \rangle ds\nonumber\\
%&-\int_0^t\big(F(\Pi_s),(\Pi_s,e_j)e_k+(\Pi_s,e_k)e_j \big)ds\nonumber\\
%&-\int_0^t\big(\sigma(\Pi_s),e_k\big)\big(\sigma(\Pi_s),e_j\big) ds
%\end{align}
%is a martingale under $\mu_0$.

For $k\in\mathbb{N}$, let $g(x)=(x,e_k)$, $x\in H$. By a similar argument, we can show that
\begin{align}\label{Mk martingale}
  & M_k(t):=g(X_t)-g(h)-\int_0^tLg(X_r)dr \nonumber\\
  =& (X_t,e_k)-(h,e_k)+\int_0^t(Ae_k,X_s)ds+\int_0^t\langle B(X_s),e_k\rangle ds-\int_0^t\big(F(X_s),e_k\big)ds
\end{align}
is a martingale under $\mu_0$.

\vskip 0.3cm
Step 3. (\ref{Mkj martingale}) and (\ref{Mk martingale}) together with It\^{o}'s formula yield that
\begin{align}\label{5.28}
<M_k, M_j>(t)=\int_0^t(\sigma(X_s),e_k)(\sigma(X_s),e_j) ds,
\end{align}
where $<M_k, M_j>$ stands for the sharp bracket of  the two martingales.
Now by  Theorem 18.12 in \cite{2002-Kallenberg-p638-638}, there exists a probability space
$(\Omega^{\prime}, {\cal F}^{\prime}, P^{\prime})$ with a filtration ${\cal F}_t^{\prime}$ such that on the standard extension
$$(\Omega\times\Omega^{\prime},  {\cal F}\times {\cal F}^{\prime}, {\cal F}_t\times {\cal F}^{\prime}_t, \mu_0\times P^{\prime} )$$
of  $(\Omega, {\cal F}, {\cal F}_t, P)$ there exists a one-dimensional Brownian motion $W(t),t\geq 0$ such that
\begin{align}\label{5.29}
M_k(t) = \int_0^t(\sigma(X_s),e_k)dW(s),
\end{align}
namely,
\begin{align}\label{5.29-1}
(X_t,e_k)-(h,e_k)=&-\int_0^t(Ae_k,X_s)ds-\int_0^t\langle B(X_s),e_k\rangle ds \nonumber\\
&+\int_0^t(F(X_s),e_k)ds+\int_0^t(\sigma(X_s),e_k)dW(s)
\end{align}
for every $k\geq 1$. Thus, under $\mu_0$, $\{X_t, t\geq 0\}$ is a solution to SNSE (\ref{Abstract}).
%\begin{align}
%X_t=h-\int_0^tAX_sds+\int_0^tb_1(X_s)ds+\int_0^tb_2(X_s)ds+\sum_{i=1}^m\int_0^t\sigma_i(X_s)dB^i_s
%\end{align}
By the uniqueness of the SNSE, we conclude that $\mu_0=\mu$ completing the proof of the theorem.$\blacksquare$

\vskip 0.3cm
In the next theorem, we will remove the restrictions placed on the coefficients and the initial value $h$.

\begin{theorem}
Assume (\hyperlink{H.1}{H.1}), (\hyperlink{H.2}{H.2}), (\hyperlink{H.3}{H.3}), (\hyperlink{H.4}{H.4}) and $h\in H$. Then,
for any $T>0$, $\mu_{\varepsilon}$ converges weakly to $\mu$, as $\varepsilon\rightarrow 0$, on the space $D([0,T], H)$ equipped with the Skorohod topology.
\end{theorem}

\noindent {\bf Proof}. For each $n\in\mathbb{N}$, let $h^n, F_n(u), \sigma_n(u), \sigma_n^{\varepsilon}(u,z)$ denote the corresponding orthogonal projections of $h, F(u), \sigma(u), \sigma^{\varepsilon}(u,z)$ into the $n$-dimensional space $\mathrm{span}\{e_1,\cdots,e_n\}$. Then, for each $n\in\mathbb{N}$, $\{\sigma_n^{\varepsilon}\}_{\varepsilon\leq\varepsilon_0}$ and $F_n$ satisfy (\hyperlink{H.2}{H.2})---(\hyperlink{H.5}{H.5}). Moreover, there is a constant $C$ independent of $n$ such that for every $u,u_1,u_2\in H$,
\begin{gather}
\label{H3 H-2}
\sup_{n\in\mathbb{N}}\|F_n(u)\|_H^2+\sup_{n\in\mathbb{N}}\|\sigma_n(u)\|_H^2 +\sup_{n\in\mathbb{N},\varepsilon\leq\varepsilon_0}\int_{\mathbb{R}_0}\| \sigma_n^{\varepsilon}(u,z)\|_{H}^2\nu(dz)\leq C(1+\|u\|_H^2),
\\
\begin{aligned}\label{H3 H Lip}
&\sup_{n\in\mathbb{N}}\|F_n(u_1)-F_n(u_2)\|_H^2+\sup_{n\in\mathbb{N}}\|\sigma_n(u_1)-\sigma_n(u_2)\|_H^2 \\ +&\sup_{n\in\mathbb{N},\varepsilon\leq\varepsilon_0}\int_{\mathbb{R}_0}\| \sigma_n^{\varepsilon}(u_1,z)- \sigma_n^{\varepsilon}(u_2,z)\|_{H}^2\nu(dz)\leq C\|u_1-u_2\|_H^2.
\end{aligned}
\end{gather}

\noindent Let $u^{n,\varepsilon}, u^n$ be the solutions of the SNSEs:
\begin{align}
\label{5.30}
u^{n,\varepsilon}(t)=& h^n-\int_0^tAu^{n,\varepsilon}(s) ds - \int_0^t B(u^{n,\varepsilon}(s)) ds +\int_0^t F_n(u^{n,\varepsilon}(s)) ds \nonumber\\
&+\int_0^t\int_{\mathbb{R}_0}\sigma_n^{\varepsilon}(u^{n,\varepsilon}(s-),z)\widetilde{N}(dzds),\\
\label{5.31}
u^n(t)=& h^n-\int_0^tAu^n(s) ds-\int_0^tB(u^n(s)) ds +\int_0^t F_n(u^n(s)) ds \nonumber\\
&+\int_0^t\sigma_n(u^n(s))dW(s).
\end{align}

\noindent By Theorem \ref{V growth theorem}, we have for each $n\in\mathbb{N}$,
\begin{align}\label{u nvare to un}
u^{n,\varepsilon}\xrightarrow{\varepsilon\rightarrow 0}u^n \quad \text{in distribution on the space $D([0,T], H)$}.
\end{align}
%In view of the continuity of $u^n$ in $H$, therefore, for each $n\in\mathbb{N}$, the finite dimensional distribution of $u^{n,\varepsilon}$ converges to the corresponding finite dimensional distribution of $u^n$.
Moreover, as the proof of (\ref{5.0}), using (\ref{H3 H-2}) we can show that
\begin{align}
\label{5.33-1}
\sup_{n\in\mathbb{N},\varepsilon\leq\varepsilon_0}\bigg\{E\sup_{0\leq t\leq T}\|u^{n,\varepsilon}(t)\|_H^4+E\Big(\int_0^T\|u^{n,\varepsilon}(s)\|_V^2ds\Big)^2\bigg\}<\infty, \\
\label{5.33-2}
\sup_{n\in\mathbb{N}}\bigg\{E\sup_{0\leq t\leq T}\|u^{n}(t)\|_H^4+E\Big(\int_0^T\|u^{n}(s)\|_V^2ds\Big)^2\bigg\}<\infty.
\end{align}

\noindent We claim  that for any $\delta>0$,
\begin{gather}
\label{5.33}
\lim_{n\rightarrow \infty}P\left(\sup_{0\leq t\leq T}\|u^{n}(t)-u(t)\|_H>\delta\right)=0,\\
\label{5.32}
\lim_{n\rightarrow \infty}\lim_{\varepsilon\rightarrow 0}P\left(\sup_{0\leq t\leq T}\|u^{n,\varepsilon}(t)-u^{\varepsilon}(t)\|_H>\delta \right)=0 .
\end{gather}
Because of similarity, we only prove (\ref{5.32}) here.
%As the proof of (\ref{5.0}), using (\ref{H3 H-2}) we can show that
%\begin{align}
%\label{5.33-1}
%\sup_{n,\varepsilon}\bigg\{E\sup_{0\leq t\leq T}\|u^{n,\varepsilon}(t)\|_H^4+E\Big(\int_0^T\|u^{n,\varepsilon}(s)\|_V^2ds\Big)^2\bigg\}<\infty, \\
%\label{5.33-2}
%\sup_{n}\bigg\{E\sup_{0\leq t\leq T}\|u^{n}(t)\|_H^4+E\Big(\int_0^T\|u^{n}(s)\|_V^2ds\Big)^2\bigg\}<\infty.
%\end{align}
Applying It\^{o}'s formula, we have
\begin{align}\label{5.34}
&e^{-\gamma \int_0^t\|u^{\varepsilon}(\rho)\|_V^2d\rho}\|u^{n,\varepsilon}(t)-u^{\varepsilon}(t)\|_H^2\nonumber\\
=& \|h^n-h\|_H^2-\gamma\int_0^t e^{-\gamma\int_0^s\|u^{\varepsilon}(\rho)\|_V^2d\rho}
\|u^{n,\varepsilon}(s)-u^{\varepsilon}(s)\|_H^2\|u^{\varepsilon}(s)\|_V^2ds\nonumber\\
& -2\int_0^te^{-\gamma\int_0^s\|u^{\varepsilon}(\rho)\|_V^2d\rho}\langle A(u^{n,\varepsilon}(s)-u^{\varepsilon}(s)),u^{n,\varepsilon}(s)-u^{\varepsilon}(s)\rangle ds\nonumber\\
&-2 \int_0^te^{-\gamma\int_0^s\|u^{\varepsilon}(\rho)\|_V^2d\rho}\langle  B(u^{n,\varepsilon}(s))-B(u^{\varepsilon}(s)),u^{n,\varepsilon}(s)-u^{\varepsilon}(s)\rangle ds\nonumber\\
&+2 \int_0^te^{-\gamma\int_0^s\|u^{\varepsilon}(\rho)\|_V^2d\rho}\big( F_n(u^{n,\varepsilon}(s))-F(u^{\varepsilon}(s)),u^{n,\varepsilon}(s)-u^{\varepsilon}(s)\big)ds\nonumber\\
&+2\int_0^t\int_{\mathbb{R}_0}e^{-\gamma\int_0^s\|u^{\varepsilon}(\rho)\|_V^2d\rho}\times \nonumber\\
&~~~~~~~~\big( \sigma_n^{\varepsilon}(u^{n,\varepsilon}(s-),z)-\sigma^{\varepsilon}(u^{\varepsilon}(s-),z), u^{n,\varepsilon}(s-)-u^{\varepsilon}(s-)\big)\widetilde{N}(dzds) \nonumber\\
&+ \int_0^t\int_{\mathbb{R}_0}e^{-\gamma \int_0^s\|u^{\varepsilon}(\rho)\|_V^2d\rho}\|\sigma_n^{\varepsilon}(u^{n,\varepsilon}(s-),z) -\sigma^{\varepsilon}(u^{\varepsilon}(s-),z)\|_H^2 N(dzds) \nonumber\\
:=& \sum_{k=1}^7 I_k^{n, \varepsilon}(t).
\end{align}
%The Ito's formula for $|X_t^{n,\varepsilon}-X_t^{\varepsilon}|_H^2$ can be found in \cite{KR}, \cite{PR}. It can also be seen through finite dimensional approximations of $X_t^{n,\varepsilon}-X_t^{\varepsilon}$.
By (\ref{b eq 0}) and (\ref{b estimate 1}) we have
\begin{align}\label{B-B simplify}
  & 2|\langle  B(u^{n,\varepsilon}(s))-B(u^{\varepsilon}(s)),u^{n,\varepsilon}(s)-u^{\varepsilon}(s)\rangle|
  = 2|\langle  B(u^{n,\varepsilon}(s)-u^{\varepsilon}(s)),u^{\varepsilon}(s)\rangle| \nonumber\\
  \leq & 4\|u^{n,\varepsilon}(s)-u^{\varepsilon}(s)\|_V \|u^{n,\varepsilon}(s)-u^{\varepsilon}(s)\|_H\|u^{\varepsilon}(s)\|_V \nonumber\\
  \leq & \|u^{n,\varepsilon}(s)-u^{\varepsilon}(s)\|_V^2 +4\|u^{\varepsilon}(s)\|_V^2\|u^{n,\varepsilon}(s)-u^{\varepsilon}(s)\|_H^2 .
\end{align}
Therefore, by (\ref{A self-adjoint}) and (\ref{B-B simplify}) we obtain that
\begin{align}\label{5.34-1}
\sum_{k=2}^4 I_k^{n, \varepsilon}(t)\leq &\int_0^te^{-\gamma \int_0^s\|u^{\varepsilon}(\rho)\|_V^2d\rho}\Big[-\|u^{n,\varepsilon}(s)-u^{\varepsilon}(s)\|_V^2 \nonumber\\ &~~~~~~+(4-\gamma)\|u^{\varepsilon}(s)\|_V^2\|u^{n,\varepsilon}(s)-u^{\varepsilon}(s)\|_H^2 \Big]ds \nonumber\\
\leq & -\int_0^te^{-\gamma \int_0^s\|u^{\varepsilon}(\rho)\|_V^2d\rho}\|u^{n,\varepsilon}(s)-u^{\varepsilon}(s)\|_V^2ds ,
\end{align}
if we take $\gamma\geq 4$. Using the Lipschitz continuity of $F$, we have
\begin{align}\label{I_5}
 &  E\sup_{0\leq s\leq t}|I_5^{n, \varepsilon}(s)| \nonumber\\
\leq & E\int_0^te^{-\gamma \int_0^s\|u^{\varepsilon}(\rho)\|_V^2d\rho}
\|u^{n,\varepsilon}(s)-u^{\varepsilon}(s)\|_H^2 ds \nonumber\\
& +E\int_0^te^{-\gamma \int_0^s\|u^{\varepsilon}(\rho)\|_V^2d\rho}
\|F_n(u^{n,\varepsilon}(s))-F(u^{\varepsilon}(s))\|_H^2 ds \nonumber\\
\leq & E\int_0^te^{-\gamma \int_0^s\|u^{\varepsilon}(\rho)\|_V^2d\rho}
\|u^{n,\varepsilon}(s)-u^{\varepsilon}(s)\|_H^2 ds \nonumber\\
& +2E\int_0^te^{-\gamma \int_0^s\|u^{\varepsilon}(\rho)\|_V^2d\rho}
\|F_n(u^{n,\varepsilon}(s))-F(u^{n,\varepsilon}(s))\|_H^2 ds \nonumber\\
& +2E\int_0^te^{-\gamma \int_0^s\|u^{\varepsilon}(\rho)\|_V^2d\rho}
\|F(u^{n,\varepsilon}(s))-F(u^{\varepsilon}(s))\|_H^2 ds \nonumber\\
\leq & CE\int_0^te^{-\gamma \int_0^s\|u^{\varepsilon}(\rho)\|_V^2d\rho}
\|u^{n,\varepsilon}(s)-u^{\varepsilon}(s)\|_H^2 ds \nonumber\\
& + 2E\int_0^te^{-\gamma \int_0^s\|u^{\varepsilon}(\rho)\|_V^2d\rho}
\|F_n(u^{n,\varepsilon}(s))-F(u^{n,\varepsilon}(s))\|_H^2 ds ,
\end{align}

\noindent By Burkh\"older's inequality, we get
\begin{align}\label{I_6}
  & E\sup_{0\leq s\leq t}|I_6^{n, \varepsilon}(s)|\nonumber\\
  \leq &  2E\Big[\int_0^t\int_{\mathbb{R}_0}e^{-2\gamma\int_0^s\|u^{\varepsilon}(\rho)\|_V^2d\rho}
\|\sigma_n^{\varepsilon}(u^{n,\varepsilon}(s),z)-\sigma^{\varepsilon}(u^{\varepsilon}(s),z)\|_H^2\times \nonumber\\ &~~~~~~~~~~~~~~~~\|u^{n,\varepsilon}(s)-u^{\varepsilon}(s)\|_H^2\nu(dz)ds\Big]^{\frac{1}{2}} \nonumber\\
\leq & 2E\bigg[\sup_{0\leq s\leq t}e^{-\frac{\gamma}{2}\int_0^s\|u^{\varepsilon}(\rho)\|_V^2d\rho}\|u^{n,\varepsilon}(s)-u^{\varepsilon}(s)\|_H \times\nonumber\\
&~~~~~~ \Big(\int_0^t\int_{\mathbb{R}_0}e^{-\gamma\int_0^s\|u^{\varepsilon}(\rho)\|_V^2d\rho}
\|\sigma_n^{\varepsilon}(u^{n,\varepsilon}(s),z)-\sigma^{\varepsilon}(u^{\varepsilon}(s),z)\|_H^2\nu(dz)ds \Big)^{\frac{1}{2}}\bigg]\nonumber\\
\leq & \frac{1}{2}E\sup_{0\leq s\leq t}e^{-\gamma\int_0^s\|u^{\varepsilon}(\rho)\|_V^2d\rho}\|u^{n,\varepsilon}(s)-u^{\varepsilon}(s)\|_H^2 \nonumber\\ &+2E\int_0^t\int_{\mathbb{R}_0}e^{-\gamma\int_0^s\|u^{\varepsilon}(\rho)\|_V^2d\rho}
\|\sigma_n^{\varepsilon}(u^{n,\varepsilon}(s),z)-\sigma^{\varepsilon}(u^{\varepsilon}(s),z)\|_H^2\nu(dz)ds \nonumber\\
\leq & \frac{1}{2}E\sup_{0\leq s\leq t}e^{-\gamma\int_0^s\|u^{\varepsilon}(\rho)\|_V^2d\rho}\|u^{n,\varepsilon}(s)-u^{\varepsilon}(s)\|_H^2 \nonumber\\
& +CE\int_0^te^{-\gamma\int_0^s\|u^{\varepsilon}(\rho)\|_V^2d\rho}\|u^{n,\varepsilon}(s)-u^{\varepsilon}(s)\|_H^2ds \nonumber\\
&+4E\int_0^t\int_{\mathbb{R}_0}e^{-\gamma\int_0^s\|u^{\varepsilon}(\rho)\|_V^2d\rho}
\|\sigma_n^{\varepsilon}(u^{n,\varepsilon}(s),z)-\sigma^{\varepsilon}(u^{n,\varepsilon}(s),z)\|_H^2\nu(dz)ds ,
\end{align}

\noindent where the uniform Lipschitz constant of $\sigma^{\varepsilon}$ has been used. Similar to (\ref{I_5}), we have
\begin{align}\label{I_7}
  & E\sup_{0\leq s\leq t}|I_7^{n, \varepsilon}(s)|\nonumber\\
  = & E\int_0^t\int_{\mathbb{R}_0}e^{-\gamma\int_0^s\|u^{\varepsilon}(\rho)\|_V^2d\rho}
\|\sigma_n^{\varepsilon}(u^{n,\varepsilon}(s),z)-\sigma^{\varepsilon}(u^{\varepsilon}(s),z)\|_H^2\nu(dz)ds \nonumber\\
\leq & CE\int_0^te^{-\gamma\int_0^s\|u^{\varepsilon}(\rho)\|_V^2d\rho}\|u^{n,\varepsilon}(s)-u^{\varepsilon}(s)\|_H^2ds \nonumber\\
&+2E\int_0^t\int_{\mathbb{R}_0}e^{-\gamma\int_0^s\|u^{\varepsilon}(\rho)\|_V^2d\rho}
\|\sigma_n^{\varepsilon}(u^{n,\varepsilon}(s),z)-\sigma^{\varepsilon}(u^{n,\varepsilon}(s),z)\|_H^2\nu(dz)ds .
\end{align}

%\begin{align}\label{sigma_nvar-sigma_var}
%& E\int_0^t\int_{\mathbb{R}_0}e^{-\gamma\int_0^s\|u^{\varepsilon}(\rho)\|_V^2d\rho}
%\|\sigma_n^{\varepsilon}(u^{n,\varepsilon}(s-),z)-\sigma^{\varepsilon}(u^{\varepsilon}(s-))\|_H^2\nu(dz)ds\nonumber\\
%\leq & 2E\int_0^t\int_{\mathbb{R}_0}e^{-\gamma\int_0^s\|u^{\varepsilon}(\rho)\|_V^2d\rho}
%\|\sigma_n^{\varepsilon}(u^{n,\varepsilon}(s-),z)-\sigma_n^{\varepsilon}(u^{\varepsilon}(s-))\|_H^2\nu(dz)ds\nonumber\\
%& +2E\int_0^t\int_{\mathbb{R}_0}e^{-\gamma\int_0^s\|u^{\varepsilon}(\rho)\|_V^2d\rho}
%\|\sigma_n^{\varepsilon}(u^{\varepsilon}(s-),z)-\sigma^{\varepsilon}(u^{\varepsilon}(s-))\|_H^2\nu(dz)ds\nonumber\\
%\leq & CE\int_0^te^{-\gamma\int_0^s\|u^{\varepsilon}(\rho)\|_V^2d\rho}\|u^{n,\varepsilon}(s)-u^{\varepsilon}(s)\|_H^2ds \nonumber\\
%&+2E\int_0^t\int_{\mathbb{R}_0}e^{-\gamma\int_0^s\|u^{\varepsilon}(\rho)\|_V^2d\rho}
%\|\sigma_n^{\varepsilon}(u^{\varepsilon}(s-),z)-\sigma^{\varepsilon}(u^{\varepsilon}(s-))\|_H^2\nu(dz)ds
%\end{align}

\noindent Combining (\ref{5.34}), (\ref{5.34-1})---(\ref{I_7}) together yields that for $t\leq T$,

\begin{align}\label{5.35}
& E\sup_{0\leq s\leq t}e^{-\gamma \int_0^s\|u^{\varepsilon}(\rho)\|_V^2d\rho}\|u^{n,\varepsilon}(s)-u^{\varepsilon}(s)\|_H^2\nonumber\\
& +2E\int_0^te^{-\gamma \int_0^s\|u^{\varepsilon}(\rho)\|_V^2d\rho}\|u^{n,\varepsilon}(s)-u^{\varepsilon}(s)\|_V^2ds\nonumber\\
\leq & 2\|h^n-h\|_H^2  + CE\int_0^te^{-\gamma \int_0^s\|u^{\varepsilon}(\rho)\|_V^2d\rho}
\|u^{n,\varepsilon}(s)-u^{\varepsilon}(s)\|_H^2 ds\nonumber\\
& + 4E\int_0^te^{-\gamma \int_0^s\|u^{\varepsilon}(\rho)\|_V^2d\rho}
\|F_n(u^{n,\varepsilon}(s))-F(u^{n,\varepsilon}(s))\|_H^2 ds \nonumber\\
& + 12E\int_0^te^{-\gamma\int_0^s\|u^{\varepsilon}(\rho)\|_V^2d\rho}\int_{\mathbb{R}_0}
\|\sigma_n^{\varepsilon}(u^{n,\varepsilon}(s),z)-\sigma^{\varepsilon}(u^{n,\varepsilon}(s),z)\|_H^2\nu(dz)ds .
\end{align}

\noindent Applying the Gronwall's inequality we obtain for $t\in [0,T]$,
 \begin{align}\label{5.36}
& E\sup_{0\leq s\leq t}e^{-\gamma \int_0^s\|u^{\varepsilon}(\rho)\|_V^2d\rho}\|u^{n,\varepsilon}(s)-u^{\varepsilon}(s)\|_H^2\nonumber\\
& +E\int_0^te^{-\gamma \int_0^s\|u^{\varepsilon}(\rho)\|_V^2d\rho}\|u^{n,\varepsilon}(s)-u^{\varepsilon}(s)\|_V^2ds\nonumber\\
\leq & C\times\Big[\|h^n-h\|_H^2 +E\int_0^t\|F_n(u^{n,\varepsilon}(s))-F(u^{n,\varepsilon}(s))\|_H^2 ds \nonumber\\
&~~~+E\int_0^t\int_{\mathbb{R}_0}
\|\sigma_n^{\varepsilon}(u^{n,\varepsilon}(s),z)-\sigma^{\varepsilon}(u^{n,\varepsilon}(s),z)\|_H^2\nu(dz)ds \Big].
\end{align}

\noindent We claim that
\begin{gather}
\label{sigma-sigma to 0}
\lim_{n\rightarrow\infty}\lim_{\varepsilon\rightarrow 0}E\int_0^T\int_{\mathbb{R}_0}
\|\sigma_n^{\varepsilon}(u^{n,\varepsilon}(s),z)-\sigma^{\varepsilon}(u^{n,\varepsilon}(s),z)\|_H^2\nu(dz)ds=0 ,\\
\label{F-F to 0}  \lim_{n\rightarrow\infty}\lim_{\varepsilon\rightarrow 0}E\int_0^T\|F_n(u^{n,\varepsilon}(s))-F(u^{n,\varepsilon}(s))\|_H^2 ds =0 .
\end{gather}

\noindent Suppose the above claims are proved. Then we conclude from (\ref{5.36}) that
\begin{align}\label{EsupH to 0}
 \lim_{n\rightarrow\infty}\lim_{\varepsilon\rightarrow 0} E\sup_{0\leq s\leq T}e^{-\gamma \int_0^s\|u^{\varepsilon}(\rho)\|_V^2d\rho}\|u^{n,\varepsilon}(s)-u^{\varepsilon}(s)\|_H^2 =0 .
\end{align}

\noindent Let us only prove (\ref{sigma-sigma to 0}). The proof of (\ref{F-F to 0}) is similar and simpler.
Let
\begin{align}
G_n^{\varepsilon}(x):=\int_{\mathbb{R}_0}
\|\sigma_n^{\varepsilon}(x,z)-\sigma^{\varepsilon}(x,z)\|_H^2\nu(dz),\quad x\in H .
\end{align}
Note that
\begin{align}\label{G_n^vare growth}
\sup_{n\in\mathbb{N},\varepsilon\leq\varepsilon_0}G_n^{\varepsilon}(x)\leq C(1+\|x\|_H^2) .
\end{align}
By (\ref{5.33-1}) and the dominated convergence theorem, to prove (\ref{sigma-sigma to 0}), it suffices to show that for each $s\in[0,T]$,
\begin{align}\label{claim2}
\lim_{n\rightarrow\infty}\lim_{\varepsilon\rightarrow 0}EG_n^{\varepsilon}(u^{n,\varepsilon}(s))=0 .
\end{align}

%We claim that for every $n\in\mathbb{N}$, $s\in[0,T]$,
%\begin{align}\label{claim2}
%\lim_{n\rightarrow\infty}\lim_{\varepsilon\rightarrow 0}EG_n^{\varepsilon}(u^{n,\varepsilon}(s))=0 .
%\end{align}
%
%\noindent By the dominated convergence theorem, it follows obviously from (\ref{G_n^vare growth}) and (\ref{claim2}) that
%\begin{align}\label{sigma-sigma to 0}
%\lim_{n\rightarrow\infty}\lim_{\varepsilon\rightarrow 0}E\int_0^T\int_{\mathbb{R}_0}
%\|\sigma_n^{\varepsilon}(u^{n,\varepsilon}(s),z)-\sigma^{\varepsilon}(u^{n,\varepsilon}(s),z)\|_H^2\nu(dz)ds=0 .
%\end{align}
%
%\noindent Using the similar methods as the proof of (\ref{sigma-sigma to 0}), we can also obatin
%\begin{align}\label{F-F to 0}
%  \lim_{n\rightarrow\infty}\lim_{\varepsilon\rightarrow 0}E\int_0^T\|F_n(u^{n,\varepsilon}(s))-F(u^{n,\varepsilon}(s))\|_H^2 ds =0 .
%\end{align}
%
%\noindent Therefore, it follows from (\ref{5.36}), (\ref{sigma-sigma to 0}) and (\ref{F-F to 0}) that
%\begin{align}\label{EsupH to 0}
% \lim_{n\rightarrow\infty}\lim_{\varepsilon\rightarrow 0} E\sup_{0\leq s\leq T}e^{-\gamma \int_0^s\|u^{\varepsilon}(\rho)\|_V^2d\rho}\|u^{n,\varepsilon}(s)-u^{\varepsilon}(s)\|_H^2 =0 .
%\end{align}

%\begin{align}
%  \lim_{n\rightarrow\infty}\lim_{\varepsilon\rightarrow 0}E\int_{\mathbb{R}_0}
%\|\sigma_n^{\varepsilon}(u^{n,\varepsilon}(s),z)-\sigma^{\varepsilon}(u^{n,\varepsilon}(s),z)\|_H^2\nu(dz)=0
%\end{align}

%Then the claim is just
%\begin{align}
%\lim_{n\rightarrow\infty}\lim_{\varepsilon\rightarrow 0}EG_n^{\varepsilon}(u^{n,\varepsilon}(s))=0
%\end{align}

\noindent Obviously, (\ref{claim2}) will follow if the following three equalities are proved.
%for every $n\in\mathbb{N}$,
\begin{align}
\label{claim2.1}
\lim_{\varepsilon\rightarrow 0}EG_n^{\varepsilon}(u^{n,\varepsilon}(s))&=\lim_{\varepsilon\rightarrow 0}EG_n^{\varepsilon}(u^{n}(s)), \quad\forall\,n\in\mathbb{N}, \\
\label{claim2.2}
\lim_{n\rightarrow\infty}\lim_{\varepsilon\rightarrow 0}EG_n^{\varepsilon}(u^{n}(s))&=\lim_{n\rightarrow\infty}\lim_{\varepsilon\rightarrow 0}EG_n^{\varepsilon}(u(s)) ,\\
\label{claim2.3}
\lim_{n\rightarrow\infty}\lim_{\varepsilon\rightarrow 0}EG_n^{\varepsilon}(u(s))&=0.
\end{align}

\noindent We first prove (\ref{claim2.1}). Since $u^n$ is a continuous process, due to (\ref{u nvare to un}), we see that for each $n\in\mathbb{N}$, $s\in[0,T]$,
\begin{align}
u^{n,\varepsilon}(s)\xrightarrow{\varepsilon\rightarrow 0} u^n(s) \quad\text{in distribution}.
\end{align}

\noindent Therefore, to prove (\ref{claim2.1}), we can use the Skorohod's representation theorem to assume that $\|u^{n,\varepsilon}(s)-u^n\|_H\rightarrow 0$ almost surely as $\varepsilon\rightarrow 0$. In view of (\ref{5.33-1}), $\{\|u^{n,\varepsilon}(s)\|_H^2\}_{\varepsilon\leq\varepsilon_0}$ is uniformly integrable, and therefore, we can further deduce that \begin{align}\label{E unvare-un to 0}
  \lim_{\varepsilon\rightarrow 0}E\|u^{n,\varepsilon}(s)-u^n(s)\|_H^2=0.
\end{align}
On the other hand,
\begin{align}\label{E Gunvare-Gun}
& E|G_n^{\varepsilon}(u^{n,\varepsilon}(s))-G_n^{\varepsilon}(u^{n}(s))| \nonumber\\
\leq & E\int_{\mathbb{R}_0}\Big|\|\sigma_n^{\varepsilon}(u^{n,\varepsilon}(s),z) -\sigma^{\varepsilon}(u^{n,\varepsilon}(s),z)\|_H^2 -\|\sigma_n^{\varepsilon}(u^{n}(s),z) -\sigma^{\varepsilon}(u^{n}(s),z)\|_H^2\Big|\nu(dz) \nonumber\\
\leq & E\int_{\mathbb{R}_0}\Big(\|\sigma_n^{\varepsilon}(u^{n,\varepsilon}(s),z) -\sigma_n^{\varepsilon}(u^{n}(s),z)\|_H +\|\sigma^{\varepsilon}(u^{n,\varepsilon}(s),z) -\sigma^{\varepsilon}(u^{n}(s),z)\|_H\Big)\nonumber\\
& \times \Big(\|\sigma_n^{\varepsilon}(u^{n,\varepsilon}(s),z) -\sigma^{\varepsilon}(u^{n,\varepsilon}(s),z)\|_H + \|\sigma_n^{\varepsilon}(u^{n}(s),z) -\sigma^{\varepsilon}(u^{n}(s),z)\|_H\Big)\nu(dz) \nonumber\\
\leq & \Big[2E\int_{\mathbb{R}_0}\Big(\|\sigma_n^{\varepsilon}(u^{n,\varepsilon}(s),z) -\sigma_n^{\varepsilon}(u^{n}(s),z)\|_H^2 \nonumber\\
&~~~~~~~~~~~~~~+\|\sigma^{\varepsilon}(u^{n,\varepsilon}(s),z) -\sigma^{\varepsilon}(u^{n}(s),z)\|_H^2\Big)\nu(dz)\Big]^{\frac{1}{2}} \nonumber\\
&\times \Big[4E\int_{\mathbb{R}_0}\Big(\|\sigma_n^{\varepsilon}(u^{n,\varepsilon}(s),z)\|_H^2 +\|\sigma^{\varepsilon}(u^{n,\varepsilon}(s),z)\|_H^2 + \|\sigma_n^{\varepsilon}(u^{n}(s),z)\|_H^2 \nonumber\\ &~~~~~~~~~~~~~~+\|\sigma^{\varepsilon}(u^{n}(s),z)\|_H^2\Big)\nu(dz)\Big]^{\frac{1}{2}} \nonumber\\
:=& I_1^{\varepsilon}\times I_2^{\varepsilon}.
\end{align}

\noindent By (\ref{H3 H-2}), (\ref{H2 H-2}), (\ref{5.33-1}) and (\ref{5.33-2}), we deduce that
\begin{align}\label{I2vare}
  \sup_{\varepsilon\leq\varepsilon_0}|I_2^{\varepsilon}|^2\leq C\sup_{n\in\mathbb{N},\varepsilon\leq\varepsilon_0}E(1+\|u^{n,\varepsilon}(s)\|_H^2+\|u^n(s)\|_H^2)< \infty.
\end{align}

\noindent (\ref{H2 H Lip}), (\ref{H3 H Lip}) and (\ref{E unvare-un to 0}) imply
\begin{align}\label{I1vare}
  |I_1^{\varepsilon}|^2\leq CE\|u^{n,\varepsilon}(s)-u^n(s)\|_H^2 \xrightarrow{\varepsilon\rightarrow 0} 0.
\end{align}

\noindent Therefore, (\ref{claim2.1}) follows from (\ref{E Gunvare-Gun}), (\ref{I2vare}) and (\ref{I1vare}). In view of (\ref{5.33}), a similar argument leads to
\begin{align}\label{claim2.2fact}
  \lim_{n\rightarrow\infty}\sup_{\varepsilon\leq\varepsilon_0}E|G_n^{\varepsilon}(u^n(s))-G_n^{\varepsilon}(u(s))|=0 .
\end{align}
Hence (\ref{claim2.2}) holds. Note that (\hyperlink{H.4}{H.4}) and the (ii) of (\hyperlink{H.3}{H.3}) imply
\begin{align}\label{sigma and sigma to 0}
  &\lim_{n\rightarrow\infty}\lim_{\varepsilon\rightarrow 0}\int_{\mathbb{R}_0} \|\sigma_n^{\varepsilon}(x,z)-\sigma^{\varepsilon}(x,z)\|_H^2\nu(dz) \nonumber\\
  = & \lim_{n\rightarrow\infty}\lim_{\varepsilon\rightarrow 0}\Big[\int_{\mathbb{R}_0}\|\sigma^{\varepsilon}(x,z)\|_H^2\nu(dz) -\int_{\mathbb{R}_0}\|\sigma_n^{\varepsilon}(x,z)\|_H^2\nu(dz)\Big] \nonumber\\
  = & \|\sigma(x)\|_H^2-\lim_{n\rightarrow\infty}\|\sigma_n(x)\|_H^2 =0 ,\quad \forall\,x\in H .
\end{align}

\noindent Therefore, (\ref{claim2.3}) immediately follows from (\ref{sigma and sigma to 0}) and (\ref{G_n^vare growth}) by the dominated convergence theorem. Thus, (\ref{sigma-sigma to 0}) is proved, and so is (\ref{EsupH to 0}).

Next, we proceed with the proof of (\ref{5.32}). For any given $\delta_1>0$, in view of (\ref{5.0}), we can choose a positive constant $M_1$ such that
\begin{align}\label{5.37-2}
&\sup_{n\in\mathbb{N},\varepsilon\leq\varepsilon_0} P\left(\sup_{0\leq t\leq T}\|u^{n,\varepsilon}(t)-u^{\varepsilon}(t)\|_H>\delta, \int_0^T\|u^{\varepsilon}(s)\|_V^2ds>M_1\right)\nonumber\\
\leq &\sup_{n\in\mathbb{N},\varepsilon\leq\varepsilon_0} P\left(\int_0^T\|u^{\varepsilon}(s)\|_V^2ds>M_1\right)\leq \delta_1.
\end{align}
On the other hand, by (\ref{EsupH to 0}), we have
\begin{align}\label{5.37-3}
&\lim_{n\rightarrow\infty}\lim_{\varepsilon\rightarrow 0} P\left(\sup_{0\leq t\leq T}\|u^{n,\varepsilon}(t)-u^{\varepsilon}(t)\|_H>\delta, \int_0^T\|u^{\varepsilon}(s)\|_V^2 ds\leq M_1\right)\nonumber\\
\leq &\lim_{n\rightarrow\infty}\lim_{\varepsilon\rightarrow 0} P\left(\sup_{0\leq s\leq T}e^{-\gamma \int_0^s\|u^{\varepsilon}(\rho)\|_V^2d\rho}\|u^{n,\varepsilon}(s)-u^{\varepsilon}(s)\|_H^2\geq e^{-\gamma M_1}\delta^2\right)\nonumber\\
\leq & e^{\gamma M_1}\frac{1}{\delta^2}\lim_{n\rightarrow\infty}\lim_{\varepsilon\rightarrow 0}E\sup_{0\leq s\leq T}e^{-\gamma \int_0^s\|u^{\varepsilon}(\rho)\|_V^2d\rho}\|u^{n,\varepsilon}(s)-u^{\varepsilon}(s)\|_H^2 = 0 .
\end{align}

\noindent Combining (\ref{5.37-2}) and (\ref{5.37-3}) together yields
\begin{align}
\lim_{n\rightarrow \infty}\lim_{\varepsilon\rightarrow 0}P\left(\sup_{0\leq t\leq T}\|u^{n,\varepsilon}(t)-u^{\varepsilon}(t)\|_H>\delta \right)\leq \delta_1.
\end{align}
Since $\delta_1$ is arbitrary, (\ref{5.32}) is proved.

\vskip 0.3cm

Finally we prove that $\mu^{\varepsilon}$ converges weakly to $\mu$. Let $\mu^{\varepsilon}_n$, $\mu_n$ denote respectively the laws of $u^{n,\varepsilon}$ and $u^{n}$ on $S:=D([0,T], H)$. Let $G$ be any given bounded, uniformly continuous function on $S$. For any $n\geq 1$,  we write
\begin{align}\label{5.38}
&\int_SG(w)\mu^{\varepsilon}(dw)-\int_SG(w)\mu(dw)\nonumber\\
=&\int_SG(w)\mu^{\varepsilon}(dw)-\int_SG(w)\mu^{\varepsilon}_n(dw)+\int_SG(w)\mu_n^{\varepsilon}(dw)-\int_SG(w)\mu_n(dw)\nonumber\\
& +\int_SG(w)\mu_n(dw)-\int_SG(w)\mu(dw)\nonumber\\
=& E[G(u^{\varepsilon})-G(u^{n,\varepsilon})]+\left( \int_SG(w)\mu_n^{\varepsilon}(dw)-\int_SG(w)\mu_n(dw)\right)\nonumber\\
& +E[G(u^{n})-G(u)].
\end{align}
Give any $\delta>0$. Since $G$ is uniformly continuous, there exists $\delta_1>0$ such that
\begin{equation}\label{5.40}
\left|E\left[G(u^{\varepsilon})-G(u^{n,\varepsilon}); \sup_{0\leq s\leq T}\|u^{n,\varepsilon}(s)-u^{\varepsilon}(s)\|_H\leq \delta_1\right]\right|\leq \frac{\delta}{4}
\end{equation}
for all $n\geq 1, \varepsilon>0$. In view of (\ref{5.32}) and (\ref{5.33}), there exists $n_1$ and then $\varepsilon_{n_1}$ such that
\begin{align}\label{5.41}
&\sup_{\varepsilon\leq\varepsilon_{n_1}}\left|E\left[G(u^{\varepsilon})-G(u^{n_1,\varepsilon}); \sup_{0\leq s\leq T}\|u^{n_1,\varepsilon}(s)-u^{\varepsilon}(s)\|_H> \delta_1\right]\right|\nonumber\\
\leq& C\sup_{\varepsilon\leq\varepsilon_{n_1}}P\left(\sup_{0\leq s\leq T}\|u^{n_1,\varepsilon}(s)-u^{\varepsilon}(s)\|_H>\delta_1\right)  \leq \frac{\delta}{4},
\end{align}
and
\begin{align}\label{5.42.1}
|E[G(u^{n_1})-G(u)]|\leq \frac{\delta}{4}.
\end{align}
On the other hand, by (\ref{u nvare to un}), there exists $\varepsilon_1>0$ such that for $\varepsilon \leq \varepsilon_1$,
\begin{align}\label{5.42}
\left|\int_SG(w)\mu_{n_1}^{\varepsilon}(dw)-\int_SG(w)\mu_{n_1}(dw)\right|\leq \frac{\delta}{4}.
\end{align}
Putting (\ref{5.38})---(\ref{5.42}) together, we obtain that for $\varepsilon \leq \min\{\varepsilon_{n_1},\varepsilon_1\}$,
\begin{align}
\left|\int_SG(w)\mu^{\varepsilon}(dw)-\int_SG(w)\mu(dw)\right|\leq \delta.
\end{align}
Since $\delta>0$ is arbitrarily small, we deduce that
\begin{align}
\lim_{\varepsilon\rightarrow 0}\int_SG(w)\mu^{\varepsilon}(dw)=\int_SG(w)\mu(dw)
\end{align}
completing the proof of the Theorem.$\blacksquare$

\vskip 0.3cm
\section{Examples}
In this section,
%we assume that $\sigma(\cdot)$ and $F(\cdot)$ satisfy (\hyperlink{H.1}{H.1}).
we give some examples of $\{\sigma^{\varepsilon}\}$ which satisfy the Hypotheses in Section 3.
\begin{proposition}\label{example of sigma}
For each $\varepsilon>0$, let
\begin{align}
\sigma^{\varepsilon}(u,z)=\sigma(\theta_{\varepsilon}(z)u)h_{\varepsilon}(z), \quad u\in H, z\in\mathbb{R}_0,
\end{align}
where $\{\theta_{\varepsilon}(\cdot)\}, \{h_{\varepsilon}(\cdot)\}$ are two families  of real-valued functions on $\mathbb{R}_0$.
Assume that $\{\theta_{\varepsilon}\}$ satisfies
\begin{align}\label{theta to 0}
\sup_{z\in\mathbb{R}_0}|\theta_{\varepsilon}(z)-1|\xrightarrow{\varepsilon\rightarrow 0} 0,
\end{align}
and $\{h_{\varepsilon}\}$ satisfies
\begin{align}
\label{h2 vare 1}  \int_{\mathbb{R}_0}|h_{\varepsilon}(z)|^2\nu(dz)\xrightarrow{\varepsilon\rightarrow 0} 1,\\
\label{h vare 0}  \sup_{z\in\mathbb{R}_0}|h_{\varepsilon}(z)|\xrightarrow{\varepsilon\rightarrow 0} 0.
\end{align}
Then $\{\sigma^{\varepsilon}\}$ satisfis (\hyperlink{H.2}{H.2})---(\hyperlink{H.4}{H.4}).
\end{proposition}

\noindent {\bf Proof}. By (\ref{theta to 0}), there exists a constant $\varepsilon_1$ such that
\begin{align}\label{theta bounded}
\sup_{\varepsilon\leq\varepsilon_1}\sup_{z\in\mathbb{R}_0}|\theta_{\varepsilon}(z)|\leq 2.
\end{align}
By (\ref{h2 vare 1}), there exists a constant $\varepsilon_2$ such that
\begin{align}\label{h2 bounded}
  \sup_{\varepsilon\leq\varepsilon_2}\int_{\mathbb{R}_0}|h_{\varepsilon}(z)|^2\nu(dz)\leq 2.
\end{align}
By (\ref{h vare 0}), there exists a constant $\varepsilon_3$ such that
\begin{align}
  \sup_{\varepsilon\leq\varepsilon_3}\sup_{z\in\mathbb{R}_0}|h_{\varepsilon}(z)|\leq 1.
\end{align}
Let $\varepsilon_0=\min\{\varepsilon_1,\varepsilon_2,\varepsilon_3\}$ and assume  $\varepsilon\leq\varepsilon_0$ in the following calculation. The linear growth condition for $\sigma$ together with (\ref{theta bounded}) and (\ref{h vare 0}) yield
\begin{align}
\sup_{\|x\|_H\leq M}\sup_{z\in\mathbb{R}_0}\|\sigma^{\varepsilon}(x,z)\|_H=&   \sup_{\|x\|_H\leq M}\sup_{z\in\mathbb{R}_0}\|\sigma(\theta_{\varepsilon}(z)x)\|_H|h_{\varepsilon}(z)| \nonumber\\
 \leq &  \sup_{\|x\|_H\leq M}\sup_{z\in\mathbb{R}_0}C(1+|\theta_{\varepsilon}(z)|\|x\|_H)\sup_{z\in\mathbb{R}_0}|h_{\varepsilon}(z)| \nonumber\\
 \leq & C(1+2M)\sup_{z\in\mathbb{R}_0}|h_{\varepsilon}(z)|\xrightarrow{\varepsilon\rightarrow 0} 0.
\end{align}
Thus, (i) of (\hyperlink{H.3}{H.3}) is satisfied. By the Lipschitz condition of $\sigma$, we have
\begin{align}\label{sigsig-sigsig to 0}
&  \bigg| \int_{\mathbb{R}_0}\big(\sigma(\theta_{\varepsilon}(z)x),e_k\big) \big(\sigma(\theta_{\varepsilon}(z)x),e_j\big)|h_{\varepsilon}(z)|^2\nu(dz) \nonumber\\
&-  \int_{\mathbb{R}_0}\big(\sigma(x),e_k\big) \big(\sigma(x),e_j\big)|h_{\varepsilon}(z)|^2\nu(dz)\bigg| \nonumber\\
\leq & \int_{\mathbb{R}_0}\big|\big(\sigma(\theta_{\varepsilon}(z)x)-\sigma(x),e_k\big) \big(\sigma(\theta_{\varepsilon}(z)x),e_j\big)\big||h_{\varepsilon}(z)|^2\nu(dz) \nonumber\\
& +  \int_{\mathbb{R}_0}\big|\big(\sigma(x),e_k\big) \big(\sigma(\theta_{\varepsilon}(z)x)-\sigma(x),e_j\big)\big||h_{\varepsilon}(z)|^2\nu(dz) \nonumber\\
\leq & \Big[\int_{\mathbb{R}_0}\|\sigma(\theta_{\varepsilon}(z)x) -\sigma(x)\|_H^2|h_{\varepsilon}(z)|^2\nu(dz)\Big]^{\frac{1}{2}} \nonumber\\
&~~~~~\times\Big[\int_{\mathbb{R}_0}\|\sigma(\theta_{\varepsilon}(z)x)\|_H^2|h_{\varepsilon}(z)|^2\nu(dz)\Big]^{\frac{1}{2}} \nonumber\\
& + \Big[\int_{\mathbb{R}_0}\|\sigma(\theta_{\varepsilon}(z)x) -\sigma(x)\|_H^2|h_{\varepsilon}(z)|^2\nu(dz)\Big]^{\frac{1}{2}} \nonumber\\
&~~~~~\times\Big[\int_{\mathbb{R}_0}\|\sigma(x)\|_H^2|h_{\varepsilon}(z)|^2\nu(dz)\Big]^{\frac{1}{2}} \nonumber\\
\leq & C\sup_{z\in\mathbb{R}_0}|\theta_{\varepsilon}(z)-1|\|x\|_H \Big[\int_{\mathbb{R}_0}|h_{\varepsilon}(z)|^2\nu(dz)\Big]^{\frac{1}{2}}
\times\bigg\{\Big[\int_{\mathbb{R}_0}(1+|\theta_{\varepsilon}(z)|^2\|x\|_H^2)\nonumber\\
& ~~~~~\times|h_{\varepsilon}(z)|^2\nu(dz)\Big]^{\frac{1}{2}} +\Big[\int_{\mathbb{R}_0}(1+\|x\|_H^2) |h_{\varepsilon}(z)|^2\nu(dz)\Big]^{\frac{1}{2}}\bigg\} \nonumber\\
& \xrightarrow{\varepsilon\rightarrow 0} 0 ,
\end{align}
where we have used (\ref{theta bounded}), (\ref{h2 bounded}) and (\ref{theta to 0}). On the other hand, (\ref{h2 vare 1}) gives
\begin{align}\label{sigsig to sigsig}
  \int_{\mathbb{R}_0}\big(\sigma(x),e_k\big) \big(\sigma(x),e_j\big)|h_{\varepsilon}(z)|^2\nu(dz)\xrightarrow{\varepsilon\rightarrow 0} \big(\sigma(x),e_k\big) \big(\sigma(x),e_j\big) .
\end{align}
Combining (\ref{sigsig-sigsig to 0}) with (\ref{sigsig to sigsig}), (ii) of (\hyperlink{H.3}{H.3}) is obtained. (\hyperlink{H.2}{H.2}) and (\hyperlink{H.4}{H.4}) can be similarly verified, we omit the details. $\blacksquare$

\begin{example} Here we give some example of $\theta_{\varepsilon}$ and $h_{\varepsilon}$.
\begin{align}
\theta_{\varepsilon}(z)=1,\  1+\varepsilon\cos z, \  1-\frac{\varepsilon}{\sqrt{2\pi}}\mathrm{e}^{-\frac{\varepsilon^2 z^2}{2}}, \cdots .
\end{align}
The following examples of $h_{\varepsilon}$ satisfy the conditions in Proposition \ref{example of sigma}.

(i) \begin{align}
h_{\varepsilon}(z)=\frac{1}{\sqrt{\nu(\{\varepsilon\leq|z|\leq 1\})}}\mathbf{1}_{\{\varepsilon\leq|z|\leq 1\}},
\end{align}
where the characteristic measure $\nu$ satisfies
\begin{align}\label{nu condition 1}
\nu(\{\varepsilon\leq|z|\leq 1\})\xrightarrow{\varepsilon\rightarrow 0}\infty, \quad\text{i.e.} \quad\nu(\mathbb{R}_0)=\infty.
\end{align}

%\begin{align}\label{nu condition 1}
%\nu(\{\varepsilon\leq|z|\leq 1\})\xrightarrow{\varepsilon\rightarrow 0}\infty.
%\end{align}

(ii) \begin{align}
h_{\varepsilon}(z)=\frac{z}{\sqrt{\int_{1\leq|z|\leq \frac{1}{\varepsilon}}|z|^2\nu(dz)}}\mathbf{1}_{\{1\leq|z|\leq \frac{1}{\varepsilon}\}},
\end{align}
where the characteristic measure $\nu$ satisfies
\begin{align}\label{nu condition 2}
\varepsilon^2\int_{1\leq|z|\leq\frac{1}{\varepsilon}}|z|^2\nu(dz)\xrightarrow{\varepsilon\rightarrow 0}\infty.
\end{align}
%For example, $\nu(dz)=\frac{1}{|z|^{1+\alpha}}dz$ for any $\alpha\in(0,2)$.

(iii) \begin{align}
h_{\varepsilon}(z)=\frac{z}{\sqrt{\int_{0<|z|\leq \varepsilon}|z|^2\nu(dz)}}\mathbf{1}_{\{0<|z|\leq \varepsilon\}},
\end{align}
where the characteristic measure $\nu$ satisfies
\begin{align}\label{nu condition 3}
\frac{1}{\varepsilon^2}\int_{0<|z|\leq \varepsilon}|z|^2\nu(dz)\xrightarrow{\varepsilon\rightarrow 0}\infty.
\end{align}

For example, if $\nu_{\alpha}(dz)=\frac{1}{|z|^{1+\alpha}}dz$, which is the characteristic measure of symmetric $\alpha$-stable processes, then for each $\alpha\in(0,2)$, $\nu_{\alpha}$ satisfies (\ref{nu condition 1}), (\ref{nu condition 2}) and (\ref{nu condition 3}).

\end{example}

\vskip 0.4cm
\noindent{\bf  Acknowledgements}\   This work is partly supported by National Natural Science Foundation of China (No.11671372, No.11431014, No.11401557).

%\noindent{\bf
%Acknowledgement}. We thank the referees for the very useful suggestions and comments.
%\vskip 0.3cm

%\bibliographystyle{unsrt}
%%\bibliographystyle{siam}

%\bibliographystyle{abbrv}
%\bibliography{20170618}

\end{document}